\documentclass[lettersize,journal]{IEEEtran}
\usepackage{amsmath} 
\usepackage{amssymb}  
\usepackage{amsfonts}
\usepackage{multirow}
\usepackage{color}
\usepackage{mathtools}
\usepackage[linesnumbered,ruled,vlined]{algorithm2e}
\usepackage{array}
\usepackage[caption=false,font=normalsize,labelfont=sf,textfont=sf]{subfig}
\usepackage{textcomp}
\usepackage{stfloats}
\usepackage{url}
\usepackage{verbatim}
\usepackage{graphicx}
\usepackage{cite}
\usepackage{float}
\usepackage{url}
\usepackage{physics}
\usepackage{stackengine}
\usepackage{tablefootnote}
\usepackage[utf8]{inputenc}
\usepackage[T1]{fontenc}
\usepackage{nomencl}

\newtheorem{remark}{Remark}

\allowdisplaybreaks

\makenomenclature
\nomenclature{\(\mathbb{N}_{\rm P},\mathbb{N}_{\rm E},\mathbb{N}_{\rm L}\)}{index set of point of common coupling (PCC) bus, EVCS bus, and non-EV load bus}
\nomenclature{\(P_k^L,Q_k^L\)}{non-EVCS load at $k^{\text{th}}$ bus}
\nomenclature{\(n,p\)}{number of bus, and number of EVCS}
\nomenclature{\(\mathcal{V},\mathcal{V}_{\rm min}\)}{Set of voltage stability index (VSI), and the minimum entry of $\mathcal{V}$}
\nomenclature{\(\mathbf{V}_k\)}{$k^{\text{th}}$ bus voltage phasor. Let $V_k^d$ and $V_k^q$ be the $d$-$q$ axis voltage.}
\nomenclature{\(P^{eD},i^{eD}\)}{demanded charging power and current vector}
\nomenclature{\(P^{e*},i^{e*}\)}{optimal charging power and current vector}
\nomenclature{\(E^{D}\)}{energy demand vector from EVs}
\nomenclature{\(\mathcal{C}^{tD},\mathcal{C}^{t\ast}\)}{charging/discharging duration at the demanded and optimal charging/discharging power rate}
\nomenclature{\(\mathcal{C}^{\mathcal{P}D},\mathcal{C}^{\mathcal{P}\ast}\)}{charging/discharging price at the demanded and optimal charging/discharging power rate}
\nomenclature{\(\mathcal{W}^{t},\mathcal{I}^e\)}{waiting time and dollar incentives}

\begin{document}

\title{Co-Optimization of EV Charging Control and Incentivization for Enhanced Power System Stability}

\author{Amit Kumer Podder, Tomonori Sadamoto, and Aranya Chakrabortty

\thanks{This work was supported in part by the US National Science Foundation under grant ECCS 1931932.}
\thanks{A. Podder and A. Chakrabortty are with the Department of Electrical and Computer Engineering, North Carolina State University, Raleigh, NC 27606. Emails: \tt\small \{apodder, achakra2\}@ncsu.edu}
\thanks{T. Sadamoto is with the Department of Mechanical Engineering and Intelligent Systems, The University of Electro-Communications, Chofugaoka, Chofu, Tokyo. Email: \tt\small sadamoto@uec.ac.jp}
}


\maketitle
\begin{abstract}
We study how high charging rate demands from electric vehicles (EVs) in a power distribution grid may collectively cause poor dynamic performance, and propose a price incentivization strategy to steer customers to settle for lesser charging rate demands so that such performance degradation can be avoided. We pose the problem as a joint optimization and optimal control formulation. The optimization determines the optimal charging setpoints for EVs to minimize the $\mathcal{H}_2$-norm of the transfer function of the grid model, while the optimal control simultaneously develops a linear quadratic regulator (LQR) based state-feedback control signal for the battery currents of those EVs to jointly improve the small-signal dynamic performance of the system states. A subsequent algorithm is developed to determine how much customers may be willing to sacrifice their intended charging rate demands in return for financial incentives. Results are derived for both unidirectional and bidirectional charging, and validated using numerical simulations of multiple EV charging stations (EVCSs) in the IEEE 33-bus power distribution model.
\end{abstract}


\printnomenclature

{\it Mathematical Notations:} We denote the sets of real and complex numbers by $\mathbb{R}$ and $\mathbb{C}$, respectively. $j:=\sqrt{-1}$. ${\bf 1}_n := [1, \ldots, 1]^{\top} \in \mathbb R^n$. Given $x_{k_1}, \ldots, x_{k_w} \in \mathbb R^{1\times N}$ and $\mathcal K := \{k_1, \ldots, k_w\}$, the quantity $[x_k]_{k \in \mathcal K} := [x_{k_1}, \ldots, x_{k_w}] \in \mathbb R^{1 \times wN}$. For $\mathbb N := \{1, \ldots, n\}$, we denote the block-diagonal matrix having matrices $M_1, \cdots, M_n$ on its diagonal blocks by ${\rm diag}(M_k)_{k \in \mathbb N}$. For $x := [x_1, \ldots, x_n]^{\top} \in \mathbb R^n$ and $\beta_1, \ldots, \beta_n \in \mathbb R$, we define $\|x\|_{2,\beta} := x^{\top}{\rm diag}(\beta_k)_{k \in \{1,\ldots,n\}}x$ and $\|x\|_{\infty,\beta} := \max\{\beta_1x_1, \ldots, \beta_nx_n\}$, where $\|x\|_{2} := x^{\top}x$. 
We define the element-wise multiplication and division by $\odot$ and $\oslash$, respectively, i.e., for $x := [x_1, \ldots, x_p]^{\top} \in \mathbb R^p$ and $y := [y_1, \ldots, y_p]^{\top} \in \mathbb R^p$, $x \odot y := [x_1y_1, \ldots, x_py_p]^{\top} \in \mathbb R^p$ and $x \oslash y := [x_1/y_1, \ldots, x_p/y_p]^{\top}$, respectively.
Given $X \geq 0$, we denote its square root $X^{\frac{1}{2}} \geq 0$ such that $X=X^{\frac{1}{2}}X^{\frac{1}{2}}$. Given a stable strictly proper transfer matrix $G(s)$, its $\mathcal H_2$-norm is defined as $\|G(s)\|_{\mathcal H_2} := \left(\frac{1}{2\pi}\int_{-\infty}^{\infty}{\rm tr}(G(j\omega)G^{\top}(j\omega))d\omega\right)^{\frac{1}{2}}$ where ${\rm tr}$ is the trace operator. For any symbol $\bullet$, we denote its setpoint as $\bullet^{*}$. For three vectors $x_1$, $x_2$, $x_3$ of same dimension, $x_1\in[x_2,x_3]$ means each element of $x_1$ is lower and upper bounded by the corresponding elements of $x_2$ and $x_3$.
\section{Introduction}\label{sec_intro}
\IEEEPARstart{A} significant amount of research has been done over the last decade on how charging of electric vehicles (EVs) should be priced cost-effectively for EV drivers, owners of EV charging stations (EVCSs), and local utility companies, considering various factors such as locations of the charging stations, traffic patterns in the neighborhoods under consideration, range anxiety, demographic factors, and locational marginal pricing of electricity during various hours of the day \cite{qian,floch,yinling,Bayram}.
Pricing includes electricity costs based on actual consumption and a portion of fixed costs for operations and maintenance of the stations such as connection fees. The price range may vary depending on the amount of energy requested by the driver, and the time that the driver allows the charging port of the EVCS to deliver that energy. Software-defined apps have been developed to estimate and display the total session fees when users enter their charging demand information. 

However, one critical issue that has been paid far less attention to is how charging controls and their associated pricing may adversely impact the dynamical characteristics and stability of the distribution grid from where the EVs are drawing power. If the number of EVs increases exponentially over the next decade, and if EV owners keep following the same charging control and pricing mechanisms as they do now, that can encourage a large fraction of drivers to charge their cars at certain specific times of the day thereby causing overloading of the grid. This, in turn, may result in poor damping of small-signal oscillations of the DC bus voltage, and voltage instability. The question is whether there are ways by which grid operators may be able to incentivize EV customers to settle for a slightly lesser charging power rate than what they desire for the same charging duration, or equivalently for a slightly longer charging duration for the same energy demand, such that these instability issues in the local distribution grid can be prevented.

This paper presents a co-optimization approach to address this open question. Stability analysis of power distribution networks with EV charging has been addressed in several recent papers. For example, the studies in \cite{eig_paper_2,eig_paper_3,new1} identify the interconnection structure between an EVCS and a DC distribution network, and other network parameters and topology as the primary factors affecting small-signal stability. Specifically, the results in \cite{new1} show how EV integration may impact small-signal oscillations in the grid. However, the analysis is limited to open-loop models only and does not address any closed-loop feedback control design. Results in \cite{eig_paper_5} introduce a comprehensive small-signal dynamic model for DC fast charging stations that encompasses both islanded mode and grid-connected mode, showing the impacts of proportional-integral (PI)-based charging controllers on dominant eigenvalues. However, the goal is to improve steady-state regulation of the charging setpoints, and not the dynamic performance of the system states. Impacts of EV integration on voltage stability have been studied in \cite{new2, new3}, but there is no clear understanding of how voltage stability can be improved via pricing and control, and what trade-offs need to be made to maintain an acceptable balance between stability margins and incentivization.

We first present a detailed dynamical model of the distribution grid with multiple EVCSs and show how the conventional PI-control-based current controllers in their power electronic circuits can easily result in poor damping responses when the EV load is high. We formulate a joint optimization and optimal control problem, where the optimization determines the optimal charging setpoints of the EVs to improve small-signal damping by minimizing the $\mathcal H_2$-norm of the grid transfer function while the optimal control part simultaneously designs a linear quadratic regulator (LQR)-based state feedback law for the battery currents to jointly minimize the risk of grid instability. A subsequent algorithm is developed to determine how much customers may be willing to sacrifice their intended charging durations in return for financial incentives. The algorithm can be implemented in the form of an app where EV drivers can submit their charging demands ahead of time to explore the possibilities of discounted rates. 
The algorithm is also extended to demonstrate the trade-off between price, small-signal stability, and voltage stability of the grid when the EVCSs operate in bidirectional charging mode. Results are validated using numerical simulations of the modified IEEE 33-bus radial distribution system model with multiple EVCSs. The benefits and drawbacks of the proposed controller on the closed-loop time response of the grid are reported. 

The rest of the paper is organized as follows. Section II recapitulates the nonlinear state-space model of the EVCS integrated distribution grid. Section III presents a motivating example to show how continued growth in EVCSs may worsen the small-signal performance and the voltage stability. Section IV presents the problem formulation and main results, first for a simplified EVCS model with unidirectional charging, and thereafter for a more generalized model with both uni- and bidirectional charging. Section V describes the practical implementation of the proposed algorithm, followed by numerical results in Section VI. Section VII concludes the paper.

 \section{Model of EVCS integrated power system}\label{sec_model}
\begin{figure}[h]
    \centering
    \includegraphics[scale=0.6]{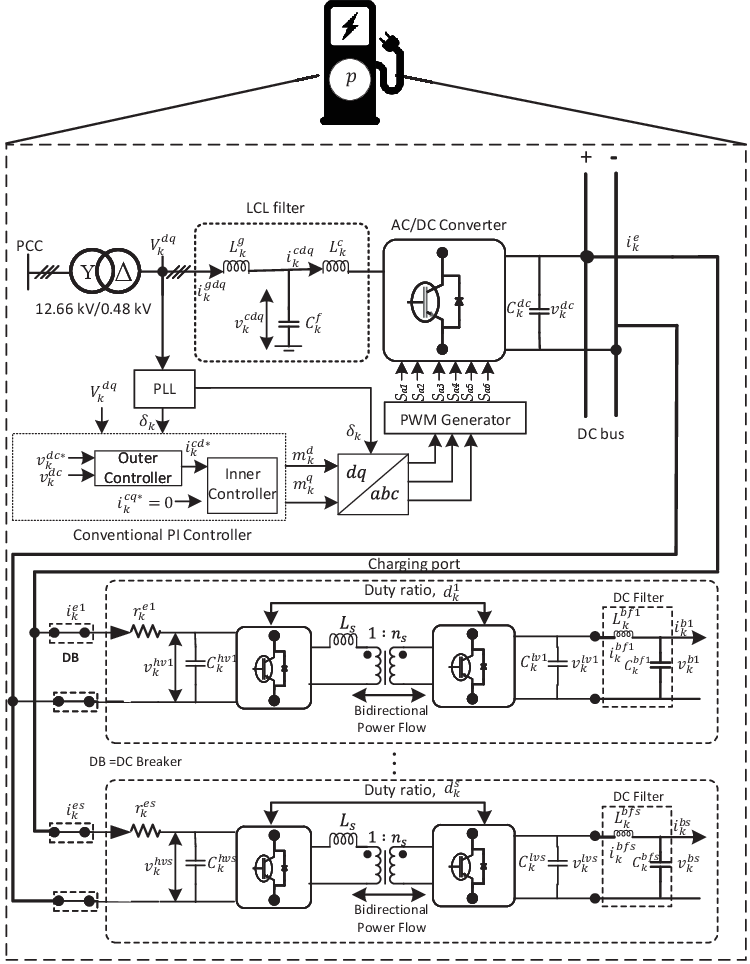}
    \caption{\footnotesize Generalized diagram of an EVCS and its internal components.}
    \label{gen_diagram}
\end{figure}

Let $\mathbb{N}$ denote the set of all buses in a radial power distribution network. We assume an aggregated generator connected at Bus 1 that supplies power to the network. This is denoted by the set $\mathbb N_{\rm P} = \{1\}$. We assume that the network consists of $p$ EVCSs, connected at buses 2 through $p+1$, denoted as $\mathbb{N}_{\rm E} = \{2, \ldots,p+1\}$. The remaining buses are assumed to be connected to non-EV loads, denoted as $\mathbb N_{\rm L} = \{p+2,\ldots n\}$. The non-EV loads are assumed to be all constant power loads that consume constant real and reactive power $P_k^L$ and $Q_k^L$, respectively, for the $k^{th}$ bus.  The assumption about the radial nature of the distribution network is made only for the simplicity of the model. This assumption does not impact our problem formulations or solution strategies.
The conventional EVCS configuration is considered, as shown in Fig.~\ref{gen_diagram},
which consists of a centralized AC/DC converter with a front-end LCL filter for creating a common DC bus, and  DC/DC converters that act as charging ports. For more details about the model, please see \cite{evcs_model, LCL_design}. 
The synchronization between the EVCS and the network is ensured by a phase-locked loop (PLL), which provides the synchronous phase angle $\delta_k$, $k\in\mathbb{N}_{\rm E}$, required for designing the controller of the converter.
\begin{table}[t]
\caption{PLL and AC/DC converter filter dynamics of $p$ EVCSs, $k\in\mathbb{N}_{\rm E}$}
\centering
\scriptsize
\begin{tabular}{|p{8cm}|}
\hline
\underline{PLL dynamics:}
\begin{equation}\label{PLL}
   \begin{cases}
    \dot{\delta}_{k}=\bar{\omega}(\omega_k-\omega_{c})\\
    \dot{\zeta}_{k}=\kappa_{k}^{I1}(V_{k}^q-0), \quad \quad 
\begin{small}
    \omega_{k}= \kappa_{k}^{P1}(V_{k}^q-0)+\zeta_{k}
\end{small}
\end{cases}
\end{equation}
where, $\delta_k$, $\omega_{k}$, $\Bar{\omega}$, and $\omega_c$  are the synchronous phase angle, angular frequency, base angular frequency, and common angular frequency, respectively. $\zeta_{k}$ is the intermediate state variable of the PLL, and $\kappa_{k}^{P1}$ and $\kappa_{k}^{I1}$ are the PI gains of the PLL. Note that $\omega_{k}$ is measured at the PLL.

\underline{AC/DC converter with LCL filter dynamics:}
\begin{equation}\label{igdq}
    \begin{cases}
        L_{k}^{g} \dot{i}_{k}^{gd}= V_{k}^{d}-v_{k}^{cd}+\omega_{k} L_{k}^{g}i_{k}^{gq}\\
        L_{k}^{g} \dot{i}_{k}^{gq}= V_{k}^{q}-v_{k}^{cq}-\omega_{k} L_{k}^{g}i_{k}^{gd}\\
    \end{cases}
\end{equation}

\begin{equation}\label{vcdq}
    \begin{cases}
      C_{k}^{f} \dot{v}_{k}^{cd}= i_{k}^{gd}-i_{k}^{cd}+\omega_{k}C_{k}^{f}v_{k}^{cq}\\
      C_{k}^{f} \dot{v}_{k}^{cq}= i_{k}^{gq}-i_{k}^{cq}-\omega_{k}C_{k}^{f}v_{k}^{cd},\\
    \end{cases}
\end{equation}

\begin{equation}\label{icdq}
    \begin{cases}
        L_{k}^{c}\dot{i}_{k}^{cd}= v_{k}^{cd}-\frac{m_{k}^{d}v_{k}^{dc}}{2}+\omega_{k} L_{k}^{c} i_{k}^{cq}\\
        L_{k}^{c}\dot{i}_{k}^{cq}= v_{k}^{cq}-\frac{m_{k}^{q}v_{k}^{dc}}{2}-\omega_{k} L_{k}^{c} i_{k}^{cd}\\
    \end{cases}
\end{equation}

where, $L_{k}^{g}$, $C_{k}^{f}$, and $L_{k}^{c}$ are the LCL filter parameters, 
$m_{k}^{d}$ and $m_{k}^{q}$ are the $d$-$q$ axis duty cycles, 
$v_{k}^{dc}$ is the DC link voltage, 
$V_{k}^{d}$ and $V_{k}^{q}$ are the $d$-$q$ axis components of the input bus voltage, 
$i_{k}^{gd}$ and $i_{k}^{gq}$ are the $d$-and $q$-axis input side filter currents, 
$v_{k}^{cd}$ and $v_{k}^{cq}$ are the voltages across the filter capacitor, and 
$i_{k}^{cd}$ and $i_{k}^{cq}$ are converter side filter currents, respectively. \\
\hline
\end{tabular}\label{filter_dynamics}
\end{table}
\begin{table}[t]
\caption{Converter controller \& DC bus dynamics of $p$ EVCSs, $k\in\mathbb{N}_{\rm E}$}
\centering
\scriptsize
\begin{tabular}{|p{8cm}|}
\hline
\underline{Outer control loop dynamics:} 
\begin{equation}\label{outercl}
\begin{cases}
\dot{\Psi}_{k}=\kappa_{k}^{I2}(v_{k}^{dc*}-v_{k}^{dc})\\ 
i_{k}^{cd*}= \kappa_{k}^{P2}(v_{k}^{dc*}-v_{k}^{dc}) +\Psi_{k}, \quad \quad i_{k}^{cq*}=0
\end{cases}
\end{equation}
\underline{Inner control loop dynamics:}
\begin{equation}\label{chid}
\begin{cases}
\dot{\chi}_{k}^{d}=\kappa_{k}^{I3}(i_{k}^{cd*}-i_{k}^{cd})\\
\begin{small}
  m_{k}^{d}=\textrm{sat}\left(\frac{2}{v_{k}^{dc}}(\kappa_{k}^{P3}(i_{k}^{cd*}-i_{k}^{cd})+\chi_{k}^{d}
+\tilde{V}_{k}^d)+\Delta m_{k}^{d}\right) 
\end{small}
\end{cases}
\end{equation}
\begin{equation}\label{chiq}
\begin{cases}
\dot{\chi}_{k}^{q}=\kappa_{k}^{I4}(i_{k}^{cq*}-i_{k}^{cq})\\
\begin{small}
 m_{k}^{q}=\textrm{sat}\left(\frac{2}{v_{k}^{dc}}(\kappa_{k}^{P4}(i_{k}^{cq*}-i_{k}^{cq})+\chi_{k}^{q}+\tilde{V}_{k}^q)+\Delta m_{k}^{q}\right),
\end{small}
\end{cases}
\end{equation}
where, $\Psi_{k}$ is the internal state variable of the outer loop, and
$\chi_{k}^{d}$ and $\chi_{k}^{q}$ are the internal state variables of the inner loop controllers, respectively. $\kappa_{k}^{P2}$, $\kappa_{k}^{I2}$, $\kappa_{k}^{P3}$, $\kappa_{k}^{I3}$, $\kappa_{k}^{P4}$, and $\kappa_{k}^{I4}$ are the  PI gains of the outer and inner loop controllers. $\tilde{V}_{k}^d=\omega_{k}L_ki_{k}^{cq}+V_{k}^d$, $\tilde{V}_{k}^q=V_{k}^q-\omega_{k} L_ki_{k}^{cd}$, and $L_k=L_{k}^{g}+L_{k}^{c}$. $\Delta m_{k}^{d}$ and $\Delta m_{k}^{q}$ are the controllable inputs, and $\textrm{sat($\cdot$)}$ is the saturation function between $[-1,1]$.\smallskip

\underline{DC bus dynamics:}
\begin{equation} \label{dc_bus}
    \dot{v}_{k}^{dc}= \frac{3(m_{k}^{d}i_{k}^{cd}+m_{k}^{q}i_{k}^{cq})}{2C_{k}^{dc}}-\frac{i_{k}^{e*}+\Delta i_{k}^{e}}{C_{k}^{dc}},
    \end{equation}
where, $C_{k}^{dc}$ is the DC link capacitance, $i_{k}^{e*}$ is the charging current setpoint, and $\Delta i_{k}^{e}$ is the controllable charging current.  \\
\hline
\end{tabular}\label{controller_dynamics}
\end{table}
The dynamics of the PLL and the AC/DC converter with LCL filter are presented in Table \ref{filter_dynamics}. The conventional PI-based controller is used for the converter, and the dynamics of the converter controller and DC bus are shown in Table \ref{controller_dynamics}. 
\begin{table}[t]
\caption{Line dynamics and power balance}
\scriptsize
\centering
\begin{tabular}{|p{8cm}|}
\hline
\underline{Line dynamics}
For $k \in \{1,\cdots,n-1\}$, and $h\in\mathcal{N}_{k}$
 \begin{equation}\label{net_1}
 \begin{cases}
         \dot{i}_{kh}^{d}=\frac{1}{l_{kh}}(V_{k}^{d}-V_{h}^{d}+\frac{G_{kh}^d\omega_cl_{kh}}{B_{kh}^d}i_{kh}^{d}+\omega_{c}l_{kh}i_{kh}^{q})\\
        \dot{i}_{kh}^{q}=\frac{1}{l_{kh}}(V_{k}^{q}-V_{h}^{q}+\frac{G_{kh}^q\omega_cl_{kh}}{B_{kh}^q}i_{kh}^{d}i_{k}^{q}-\omega_{c}l_{kh}i_{k}^{d}), 
 \end{cases}
 \end{equation}
where, \{$i_{kh}^{d}$,$i_{kh}^{q}$\} and \{$V_h^{d}$, $V_h^{q}$\} are  
the $d$-$q$ axis branch currents and bus voltages neighboring to the $k^{\rm th}$ bus, respectively. $l_{kh}$ is the line reactance.
\\\\
\underline{Power balance}
At PCC, i.e., for $k\in\mathbb{N}_{\rm P}$,
 {\scriptsize
  \begin{equation}\label{gen}
  \hspace{-3mm} \begin{cases}
  P^g_{k}=\frac{3}{2}\left(\sum_{h \in \mathcal N_k}V_{k}^dV_{h}^d G_{kh}^d+\sum_{h \in \mathcal N_k} V_{k}^qV_{h}^q G_{kh}^q\right)
  \\
  \begin{small}
  Q^g_{k}=\frac{3}{2}\left(\sum_{h \in \mathcal N_k} V_{k}^q V_{h}^d B_{kh}^d-\sum_{h \in \mathcal N_k}V_{k}^d V_{h}^q B_{kh}^q\right).
  \end{small}
  \end{cases}
        \end{equation}}
For $k\in\mathbb{N}_{\rm L},\mathbb{N}_{\rm E}$,
 {\scriptsize
  \begin{equation}\label{ev-nonev}
  \begin{cases}
   P^{L}_{k}+P_{k}^{e}=\frac{3}{2}\left(\sum_{h \in \mathcal N_k}V_{k}^dV_{h}^d G_{kh}^d+\sum_{h \in \mathcal N_k} V_{k}^qV_{h}^q G_{kh}^q\right)
  \\
  \begin{small}
  Q_{k}^{L}+Q_{k}^{e}=\frac{3}{2}\left(\sum_{h \in \mathcal N_k} V_{k}^q V_{h}^d B_{kh}^d-\sum_{h \in \mathcal N_k}V_{k}^d V_{h}^q B_{kh}^q\right)
  \end{small}
  \end{cases}
        \end{equation}}
where, 
\begin{equation}\label{PQke}
    P_{k}^{e}=\frac{3}{2} (V_{k}^di_{k}^{gd}+V_{k}^qi_{k}^{gq}),~
   Q_{k}^{e}=\frac{3}{2} (V_{k}^{q}i_{k}^{gd}-V_{k}^di_{k}^{gq}) 
\end{equation}
$P_{k}^{g}$ and $Q_{k}^{g}$ in \eqref{gen} are the active and reactive power generated by the bulk transmission grid (modeled as a constant voltage source), $P_{k}^{L}$ and $Q_{k}^{L}$ in \eqref{ev-nonev} are that consumed by the non-EV loads, while $P_{k}^{e}$ and $Q_{k}^{e}$ in \eqref{PQke} are that consumed by the EVs.\\
\hline
\end{tabular}\label{net_power_balance}
\end{table}
Each EVCS may have multiple charging ports connected to their respective DC buses, where every charging port has similar internal dynamics \cite{eig_paper_5}. The dynamics of the charging ports are ignored for simplicity as their time scale is significantly faster than that of the grid state variables that we are interested in. This will also be verified later via an eigenvalue analysis in the motivating example of Section \ref{sec_motiv}. 
Short-line models are considered for the interconnecting tie-lines between the buses. The $d$-$q$ axis admittance of the line connecting the $k^{th}$ bus and the $h^{th}$ bus, $k=1,\dots, n-1$, is denoted as $G_{kh}^{dq}+jB_{kh}^{dq}\in \mathbb{C}^2$, for any $h \in \mathcal N_k$, where $\mathcal N_k$ is the set of neighboring buses for the $k^{th}$ bus, the conductance vector $G_{kh}^{dq}:=[G_{kh}^{d}, G_{kh}^{q}]^{\top}$, and the susceptance vector $B_{kh}^{dq}:=[B_{kh}^{d}, B_{kh}^{q}]^{\top}$. The line dynamics and power flow using $d$-$q$ axis components are presented in Table \ref{net_power_balance}.

The overall state-space model of the EVCS-integrated power system is presented in a compact form by the following differential-algebraic equations (DAEs):  
\begin{equation}\label{overall_ss}
   \Sigma:\begin{cases}
   \dot{x}=f(x,y,u,\alpha)\\
   0=g(x,y,\alpha),\\
\end{cases}
\end{equation}
where, from Tables \ref{filter_dynamics}-\ref{net_power_balance}, 
\begin{equation} \label{eq1}\nonumber
\begin{split}
 x &:=\hspace{-1mm} [[(x_{k}^e)^\top]_{k \in \mathbb{N}_{\rm E}}^{\top},\hspace{-1mm}[i_{kh}^{d},i_{kh}^{q}]^{\top}_{k \in \{1,\ldots,n-1\},h\in\mathcal{N}_k}]^{\top}\hspace{-2mm} \in \hspace{-1mm}\mathbb R^{12p+2(n-1)}\\
 x_{k}^{e}  &:=\hspace{-0.5mm} \scalebox{0.98}{$ [\delta_{k},\zeta_{k}, i_{k}^{gd}, i_{k}^{gq}, i_{k}^{cd}, i_{k}^{cq}, v_{k}^{cd}, v_{k}^{cq} ,\Psi_{k},\chi_{k}^{d},\chi_{k}^{q}, v_{k}^{dc}]^{\top}\hspace{-1mm} \in\mathbb{R}^{12} $}\\
u&:=\hspace{-0.5mm}\left[\Delta m_{k}^{d},\Delta m_{k}^{q},\Delta i_{k}^e\right]^{\top}_{k \in \mathbb N_{\rm E}}\in\mathbb{R}^{3p}\\
y&:=\hspace{-0.5mm}\left[V_{k}^{d},V_{k}^{q}\right]^{\top}_{k \in \{1,\ldots,n\}}\in\mathbb{R}^{2n}, \;\alpha:=[i_k^{e*}]^{\top}_{k\in\mathbb{N}_{\rm E}}\in\mathbb{R}^p.
\end{split}
\end{equation}
The functions $f(\cdot)$ and $g(\cdot)$ follow from (\ref{PLL}) to (\ref{ev-nonev}). The equilibrium of this model is determined by substituting $f(\cdot)=g(\cdot)=0$, and solving for ($x^*,\,y^*,\,u^*$). The algebraic variable $y$ is eliminated and the model is then linearized around the equilibrium point $(x^*,\,u^*)$ by Jacobian linearization to obtain a small-signal state-variable model that will be used shortly in the next section. Note that the control input $u$ is already in the small-signal form, i.e., $u^* = 0$.

\begin{figure}[h]
    \centering
     \includegraphics[width=0.35\textwidth]{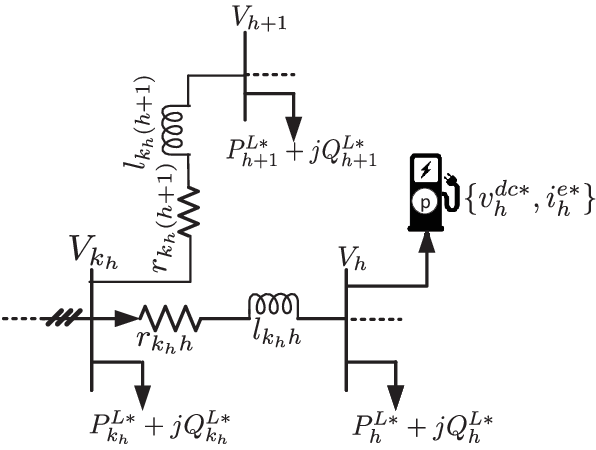}
    \caption{ \footnotesize A schematic of the radial network topology showing sending and receiving end buses.}
    \label{network}
\end{figure}

We close this section by recalling the definition of voltage stability index (VSI) that will be used shortly in the design reported in Section \ref{section_iii_C}. The VSI represents the voltage stability margin at any given bus of the distribution network and is defined at each receiving end bus. For a radial network, every receiving end bus only has one sending end bus, as shown in Fig. \ref{network} (please note that the reverse may not be true, as the sending end bus may provide power to multiple receiving end buses). 
For $h \in \{2, \ldots, n\}$ (which is the index of the receiving end bus), let the sending bus index be denoted as $k_h \in \{1, \ldots, n-1\}$. Following \cite{ref1}, we define the VSI at the $h^{\text{th}}$ bus as
\begin{align}
    &\scalebox{1.0}{$  \mathcal{V}_{h} := V_{k_h}^{*4}-4((P^{L*}_{h}+ v^{dc*}_{h}i_{h}^{e*})r_{k_hh}-Q^{L*}_{h}\mathcal{X}_{k_hh})^2 $} \nonumber \\
    &\scalebox{1.0}{$ -4((P^{L*}_{h}+ v^{dc*}_{h}i_{h}^{e*})r_{k_hh}+Q^{L*}_{h}\mathcal{X}_{k_hh})V_{h}^{*2},$} \label{vsi_main_exp}
\end{align}
where, $V_{k_h}$ is the voltage at the corresponding sending end bus, $r_{k_hh}$ and $\mathcal{X}_{k_hh}:= 2\pi\bar{\omega}l_{k_hh}$ are the resistance and reactance between $k_h^{\text{th}}$ and $h^{\text{th}}$ buses, respectively. $P^{L*}_{h}$, and $Q^{L*}_{h}$ are steady-state real and reactive powers consumed by the non-EV load in the network at the $h^{\text{th}}$ bus. Note that $\mathcal V_{h} \in [0,1]$ from the definition. The closer the VSI is to unity, the less susceptible that bus is to a voltage collapse. 
\section{Motivating Example}\label{sec_motiv}
\begin{figure}[h]
    \centering
     \includegraphics[scale=0.45]{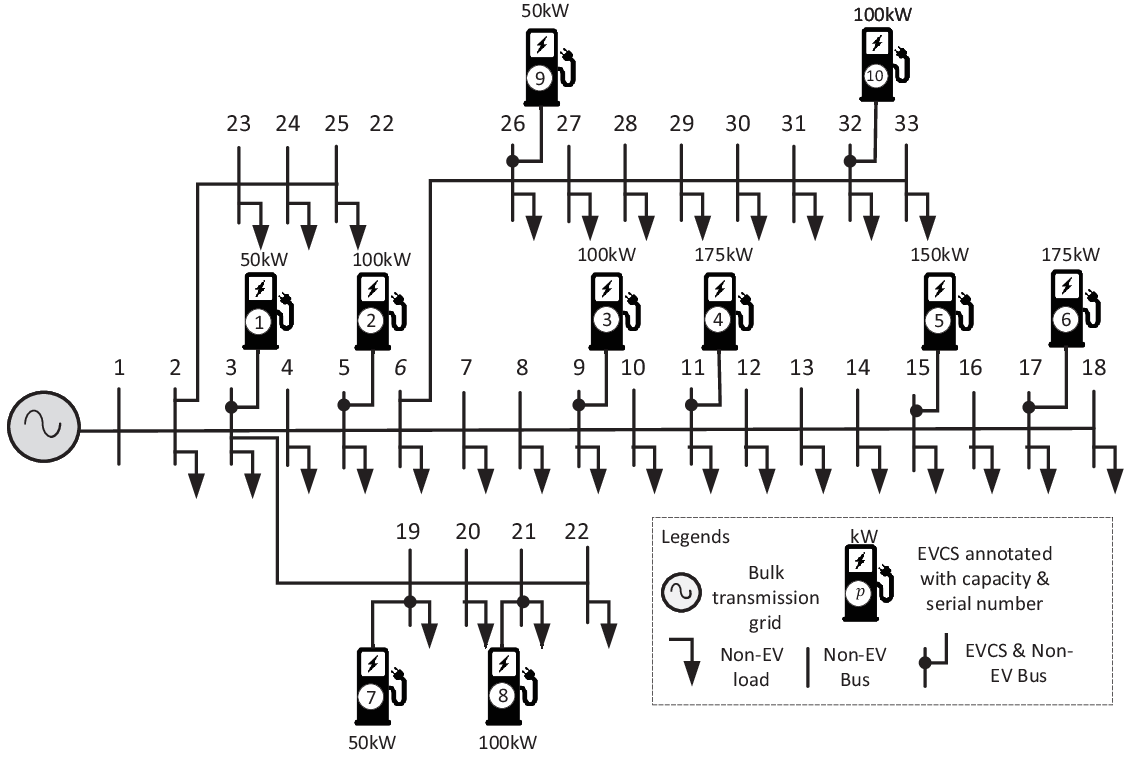}
    \caption{ \footnotesize  IEEE 33-bus radial distribution system model with $p=10$ EVCSs. }
    \label{33_bus}
\end{figure}

Consider the IEEE 33-bus medium voltage power distribution network model, where $p$ EVCSs are connected to the network as shown in Fig.~\ref{33_bus}.
Three cases are considered to evaluate the impact of increasing the number of EVCSs in this model, namely, $p=3,\,5,\,10$, with the respective EVCS-connected buses represented as $\mathbb{N}_{\rm E}:=\{3,19,5\}$, $\mathbb{N}_{\rm E}:=\{3,5,9,19,21\}$, $\mathbb{N}_{\rm E}:=\{3,5,9,11,15,17,19,21,26,32\}$. In all cases, the EVCS capacity is assumed to be between 50 kW and 175 kW. For this motivating example, we define $\Sigma_{\rm detail}$\footnote{$\Sigma_{\rm detail}$ is used only for the motivating example in this section so that we can show the entire range of eigenvalues that one may encounter in a realistic power system. For our algorithm design in Section \ref{sec_proposed_alg}, we will use the simplified model \eqref{overall_ss} that does not include the charging port dynamics.} as the system that includes not only the dynamics of \eqref{overall_ss} but also the dynamics of the dual-active bridge (DAB)-based DC/DC converters that act as the charging ports \cite{eig_paper_5}. The model parameters for a typical 50 kW EVCS are enlisted in the Appendix. Other relevant model parameters are considered following \cite{LCL_design}.
\begin{figure*}[ht]
\centering
\includegraphics[width=0.95\textwidth]{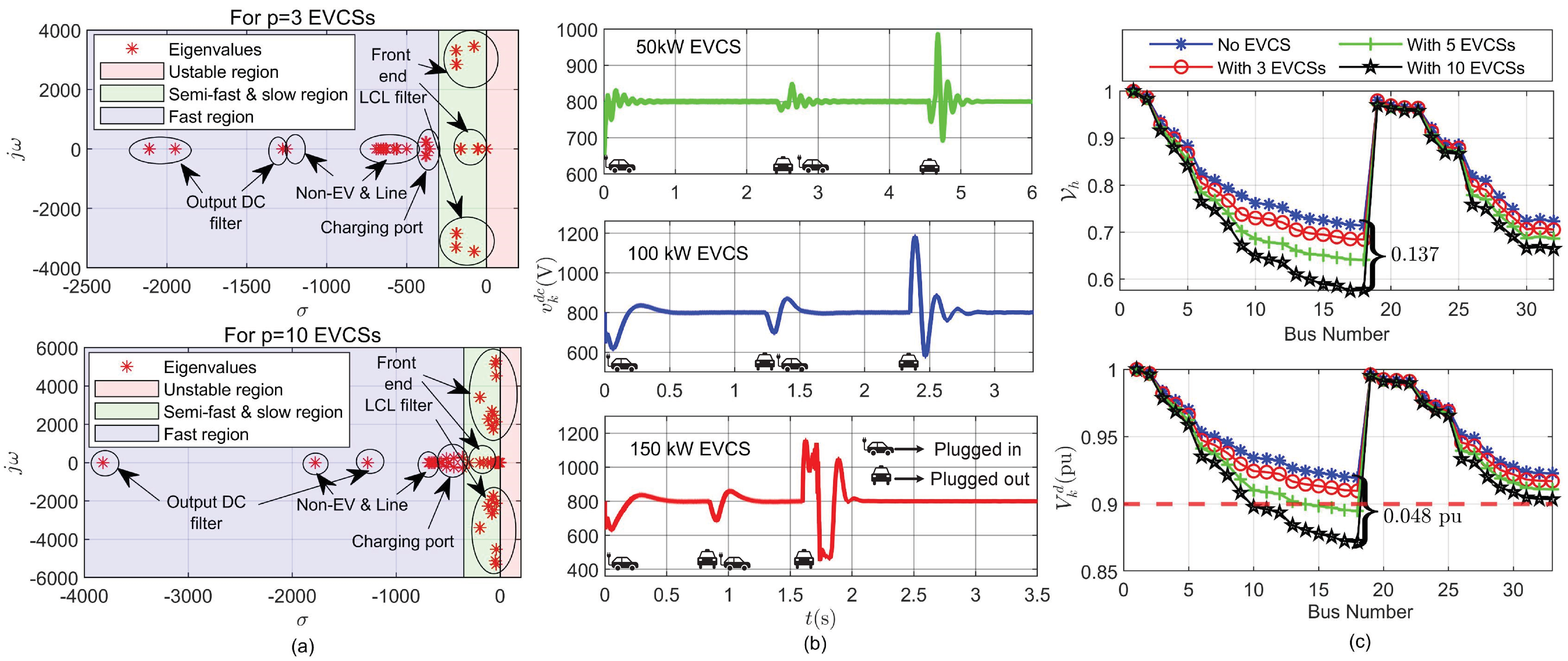}
\caption{\footnotesize (a) Open-loop eigenvalues of the linearized 33-bus distribution grid model $\Sigma_{\rm detail}$ for $p=3$ and $10$ EVCSs. (b) The transient DC bus voltage $v_k^{dc}$ inside the $1^{\text{st}}$, $3^{\text{rd}}$, and $5^{\text{th}}$ EVCSs for $p=10$, respectively. Note that the charging duration is scaled down from the real-time charging time to highlight the transients in a shorter duration. (c) Trends of VSI values of different buses of the network with the increasing number of EVCSs. The red dotted line shows the minimum allowable voltage limit in the network.}
\label{eig_time_vsi}
\end{figure*}
Figure~\ref{eig_time_vsi}(a) shows the eigenvalues of the linearized model of $\Sigma_{\rm detail}$ for $p=3$ and $p=10$ EVCSs. It can be seen that the eigenvalues with the highest participation factor from the charging port states, the output DC filter states, and the line states (all marked in the purple region) have much faster time constants than those in the green region. We label the purple eigenvalues as ``fast'', and the green eigenvalues as ``dominant''. The dominant eigenvalues show poor damping factors. The damping worsens as the number of EVCSs increases. The transient responses of the DC bus voltage $v_{k}^{dc}$ for a selected set of EVCSs, triggered by the plugging-in and out of their EV loads are shown in Fig.~\ref{eig_time_vsi}(b). The figure shows that the voltage suffers from continued oscillations that become more severe with the increase in EV charging rates and demands. Similarly, the VSI values at the different buses are shown in Fig.~\ref{eig_time_vsi}(c). It can be seen that the VSI decreases with increasing EV demand. For example, for  $p=10$, the minimum VSI and bus voltage at Bus 18 drop-down by 0.137 pu and 0.048 pu, respectively, compared to their values for $p=3$. In summary, increasing the number of EVCSs impacts both small-signal and voltage stability adversely. 

In the following sections, we resolve this problem in two ways. The first approach is to request the EV customers to reduce their demanded charging rates in exchange for an incentive. The challenge of achieving a balance between these incentives and the small-signal performance of the grid is discussed in Sections \ref{section_iii_A} and \ref{section_iii_B}. The second approach is a more comprehensive solution that takes into account not only incentivization but also leverages the bidirectional charging capability of the EV batteries. This is presented in sections \ref{section_iii_C} and \ref{section_iii_D}. 

\section{Proposed Algorithms}\label{sec_proposed_alg}
\subsection{Problem 1: Simplified EVCS integrated grid}\label{section_iii_A}
We assume that a central authority, referred to here onward as the central point operator (CPO), is in charge of running all the EVCSs in the grid. Let us assume that on any given day, $\mathcal{N}_v$ number of cars submit information about what time of the day they are planning to arrive at the $k^{th}$ EVCS to charge their EVs the following day, and how much Kilo-Watt-Hour (kWh) of energy they need for charging. 
Let the energy demand vector be denoted as $E^D:=[E_k^D]_{k \in \mathbb{N}_{\rm E}}^{\top}\in\mathbb{R}^p$, from which the EVCS owners compute the demanded charging current rate vector $i^{eD}:=[i_k^{eD}]_{k \in \mathbb{N}_{\rm E}}^{\top}\in\mathbb{R}^p$. We assume that each EVCS has a fixed capacity and only operates in grid-to-vehicle (G2V) mode, which implies that the total amount of steady-state current that the EVCSs can draw from the grid will always remain constant, irrespective of which car is coming to which port and when. 
In other words, for a given set of vehicles, $i^{eD}$ is always equal to a constant number, regardless of the schedule of charging these vehicles. The demands depend on the battery specification of each specific vehicle as well as other factors depending on traffic patterns of neighborhoods, constituting the so-called {\it customer's objective} \cite{yinling}. Our goal of this paper, in contrast, is to consider the {\it CPO's objective}, which is to determine the ideal value of the setpoint vector $i^{e*}:=[i_{k}^{e*}]_{k \in \mathbb N_{\rm E}}^{\top}\in\mathbb{R}^p$ (which may be less than the corresponding demanded value $i^{eD}$) as well as to design a state-feedback controller for $u(t)$ in \eqref{overall_ss} to maximize the closed-loop damping performance of the EVCS-integrated grid against any disturbance. From the physics of power system models, it is intuitive that the grid will always have a better small-signal stability margin when $i^{e*}$ is small versus high, i.e., when the EV loading is minimal. This indicates that the solution to the damping maximization problem will always tend towards the allowable lower bound for the demand. To avoid this problem, we add an {\it incentivization} term to the optimization objective that will allow customers to settle for a demand that is higher than the lowest allowable demand but lower than their desired demand $i^{eD}$ in exchange for money. The practical implication is that the grid operator must incentivize the customers for the deficit amperes that they are sacrificing for the sake of enhancing the stability of the grid. The problem is formulated as: 
\begin{align}
\min_{K,i^{e*}}J= & (1-\gamma) \underbrace{\int_{0}^{\infty} \|\hat{z}\|_2^2 dt}_{J_1} +\gamma \underbrace{\|i^{e*}-i^{eD}\|_{2,\beta}^2}_{J_2} \label{main_problem} \\
\text{s.t} \quad & (\ref{overall_ss}) \;\mbox{holds}, \\
&\hat{z}:= \mathcal{Q}\hat{x} + \mathcal{R}u\label{LQR}, \quad u=-K\hat{x}, \\
& i_k^{e*} \in [\max\{\underbar{$i$}_k^{e*}, 0\}, i_k^{eD}], \quad k \in \mathbb N_{\rm E}, \label{ie_star}
\end{align}
where, $\hat{x}:=P(x-x^{*}) \in \mathbb R^{7p}$ is a reduced-dimensional state vector that consists of only the physical state variables of the EVCSs, and excludes the five non-physical states per EVCS, namely,  $(\delta_{k},\zeta_{k})$ coming from the PLL, and $(\Psi_{k},\,\chi_{k}^{d},\,\chi_{k}^{q})$ coming from the internal PI controllers of the converter, as well as the two line states $(i_{kh}^{d},\,i_{kh}^{q})$ per line. This dimensionality reduction is motivated by the observations made from numerical studies, including the results presented in Section III, that show that the PLL states, controller states, and line states have the lowest participation factor in the dominant eigenvalues of the small-signal model of \eqref{overall_ss}. For example, for the IEEE 33-bus model, these states are consistently found to have less than 15\% participation in the dominant eigenvalues for $p=3,\,5,$ and $10$. An equivalent interpretation is that the time constants associated with these states are notably smaller than those with the physical states of the EVCSs, indicating that they can be eliminated using a time-scale separation from the latter. Elimination of these unimportant states facilitates both computation and communication needed for the state feedback. The matrix $P \in \mathbb R^{7p \times (12p+2(n-1))}$ is constructed accordingly using singular perturbation theory \cite{khalil_book}, resulting in $\hat{x}:=[i_{k}^{gd},\, i_{k}^{gq}, \,i_{k}^{cd}, \,i_{k}^{cq}, \,v_{k}^{cd},\, v_{k}^{cq}, \, v_k^{dc}]^{\top}$ for $k\in \mathbb N_{\rm E}$. We define
$\mathcal{Q}:=\begin{bmatrix} Q^{\frac{1}{2}} & 0 \end{bmatrix}^{\top}$, $\mathcal{R}:=\begin{bmatrix} 0 & R^{\frac{1}{2}} \end{bmatrix}^{\top}$, where $Q \succeq 0\in\mathbb{R}^{7p\times7p}$ and $R \succ 0\in\mathbb{R}^{3p\times3p}$ are constant weight matrices that can be decided by the grid operator depending on how much damping improvement is needed versus how much control energy can be spent on the feedback.  $\beta_k$ is the price incentive constant,  $\underbar{$i$}_{k}^{e*}$ is the $k$-th lowest allowable demand, and $\hat{z}$ represents the measure for quantifying the damping performance. As all EVCSs are assumed to be unidirectional, the max function is used in  \eqref{ie_star}. The matrix $K\in\mathbb{R}^{3p\times7p}$ is the optimal control gain, and $\gamma \in [0,\,1]$ is a weighting factor. 
 
\subsection{Solution Approach of Problem 1}\label{section_iii_B}
 For a given $i^{e*}$, a reduced-order linearized model of \eqref{overall_ss} with the performance output $\hat{z}$ is uniquely determined as
\begin{equation}\label{h2pfm_ss_original}
    G_{i^{e*}}:~\begin{cases}
         \dot{\hat{x}}=A(i^{e*})\hat{x} + B(i^{e*})u, \;\; \hat{x}(0)=\hat{x}_0\\
        \hat{z} = \mathcal Q\hat{x} + \mathcal R u,
    \end{cases}
\end{equation}
where $\mathcal Q$, $\mathcal R$ are defined as in \eqref{LQR}. A relaxation strategy is used to make $K$ uniquely determined as a function of $i^{e*}$ as:
\begin{equation}
   \min_{i^{e*}} J_{R}=  (1-\gamma)\underbrace{ \|G_{i^{e*}}\|_{\mathcal{H}_2}^2}_{J_{R1}} + \gamma\underbrace{\|i^{e*}-i^{eD}\|_{2,\beta}^2}_{J_{R2}} \label{relax}
\end{equation}
\begin{align}
\text{s.t} \quad& \eqref{ie_star}, (\ref{h2pfm_ss_original}) \;\mbox{hold},\\
& u = -K(i^{e*})\hat{x}, \label{Klin}\\
&K(i^{e*})={\rm LQR}(i^{e*}) \label{Klin_val},
\end{align}\text{where},

 \begin{align}\label{relax3}
     {\rm LQR}(i^{e*}) &={\underset{K(i^{e*})}{\arg \min}}\norm{\begin{bmatrix}
         Q^{\frac{1}{2}} \\  R^{\frac{1}{2}} 
     \end{bmatrix} (sI-A_{cl}(i^{e*}))^{-1}}_{\mathcal{H}_2}^{2}, \\
     A_{cl}(i^{e\ast}) &:= A(i^{e\ast})-B(i^{e\ast})K(i^{e*})\in{\mathbb{R}^{7p\times7p}}. \label{defAcl}
 \end{align} 
To find the balance of the optimal setpoint $i^{e*} \in \mathbb R^p$ considering the damping objective ($J_{R1}$) and the incentivization objective ($J_{R2}$), we propose a gradient-based algorithm by using the concept of minimizing the $\mathcal{H}_2$-norm of a transfer function as presented in \cite{h2pfm} with a minor modification. The procedure is summarized as follows.

The transfer function matrix $G_{i^{e\ast}}(s)$ in (\ref{h2pfm_ss_original}) with the controller in \eqref{Klin} for the initial state disturbance $\hat{x}(0)$ can be written as
\begin{equation}\label{transfer_function}
    G_{i^{e\ast}}(s) := C(sI-A_{cl}(i^{e*}))^{-1}\hat{x}_0,
\end{equation}
where, $C := \mathcal Q + \mathcal R K$. Let 
$\ell$ be the iteration number and $i^{e*(\ell)}$ be the $\ell$-th guess for an optimal $i^{e*}$. 
Let $p_k$ be the index corresponding to the EVCS connected to the $k$-th bus (e.g., if $\mathbb{N}_{\rm E}:=\{3,19,5\}$, then $p_3 = 1$, $p_{19} = 2$). Under this setting, we first consider deriving a linear perturbed model of $G_{i^{e*}}$ for the decision variable $i^{e*}$, described as follows. 
Given $\epsilon > 0$, for $k \in \mathbb N_{\rm E}$, we consider a perturbed setpoint by increasing only the $p_k$-th element of $i^{e*(\ell)}$ as 
\begin{align}\label{cl_state_1}
    i_k^{e*(\ell)} := i^{e*(\ell)} + e_k^p \epsilon, \quad e_k^p := [0, \ldots, 1, \ldots, 0]^{\top} \in \mathbb R^p. 
\end{align}
Then, the closed-loop state matrix at this perturbed setpoint can be written as 
\begin{equation}\label{cl_state_2}
   A_{cl,k}:=A(i_k^{e*(\ell)})-B(i_k^{e*(\ell)})K(i_k^{e*(\ell)}) 
\end{equation}
 where $K$ is designed from \eqref{Klin_val}-\eqref{relax3}. Therefore, by denoting the difference of the transfer function as
\begin{align}\label{defG_delta}
    \Delta G_k := \left(C(sI-A_{cl,k})^{-1}\hat{x}_{0,k} - G_{i^{e*(\ell)}}\right) / \epsilon, 
\end{align}
the linear approximation of $G_{i^{e\ast}}$ around $i^{e*} = i^{e*(\ell)}$ can be written as 
\begin{align}\label{def_app}
 G_{i^{e\ast}}(s) \approx G_{i^{e*(\ell)}}(s) + \textstyle \sum_{k \in \mathbb N_{\rm E}}\Delta G_{k}(s) \epsilon. 
\end{align}
Define $\{A_{\Delta, k}, B_{\Delta, k}, C_{\Delta, k}\}$ such that $\Delta G_{k} = C_{\Delta, k}(sI - A_{\Delta, k})^{-1}B_{\Delta, k}$. Then, the RHS of \eqref{def_app} can be written as $\dot{\mathcal X} = \mathcal A \mathcal X + (\mathcal B_1 + \mathcal B_2 \epsilon {\bf 1}_p)\varrho$, $\hat{z} = \mathcal C \mathcal X$ where $\varrho$ is the Dirac delta function, $\mathcal{B}:=\mathcal B_1 + \mathcal B_2 \epsilon {\bf 1}_p$, and 
\begin{align}
 &\mathcal A := {\rm diag}(A_{cl}(i^{e*(\ell)}), {\rm diag}(A_{\Delta,k})_{k \in \mathbb N_{\rm E}}), \label{defAAA}\\
 &\mathcal B_1 := [\hat{x}_0^{\top}, 0]^{\top}, \\
 &\mathcal B_2 := [0, {\rm diag}(B_{\Delta,k}^{\top})_{k \in \mathbb N_{\rm E}}]^{\top}, \\
 &\mathcal C := [C, [C_{\Delta, k}]_{k \in \mathbb N_{\rm E}}]. \label{defCCC}
\end{align}
Following \cite{h2pfm}, the gradient of the objective function $J_R$ in (\ref{relax}) is then determined as 
\begin{equation}\label{grad_obj}
    \left.\frac{\partial J_{R}}{\partial i^{e\ast}}\right|_{i^{e\ast}=i^{e\ast(\ell)}}= (1-\gamma)2\mathcal B_2^{\top}L\mathcal B_1+2\gamma\mathcal D (i^{e\ast(\ell)}-i^{eD}), 
\end{equation}
where, $\mathcal D := {\rm diag}(\beta_k)_{k \in \mathbb N_{\rm E}}$ and $L \geq 0$ is the solution of the Lyapunov equation
\begin{equation}\label{lyp}
    L\mathcal A+\mathcal A^{\top}L+\mathcal C^{\top}\mathcal C=0.
\end{equation}
Detailed derivations are given in the Appendix. After obtaining the gradient, the parameter of interest $i^{e\ast}$ is updated as 
\begin{equation}\label{updateeq1}
    i^{e\ast(\ell+1)}=i^{e\ast(\ell)} - \alpha_s \left.\frac{\partial J_{R}}{\partial i^{e\ast}}\right|_{i^{e\ast}=i^{e\ast(\ell)}},
\end{equation}
where, $\alpha_s > 0$ is the learning rate of the gradient. To find the optimal value of $\alpha_s$ the traditional Armijo-rule-based line search is utilized which guarantees finite iteration convergence of (\ref{relax})-(\ref{Klin_val}) to a local optimum or the boundary of the constraint set. The steps for solving (\ref{relax})-(\ref{Klin_val}) are summarized in Algorithm \ref{alg_1}.

\begin{algorithm}[h]
\caption{Algorithm for joint parametric optimization and optimal control (\ref{relax})-(\ref{Klin_val})}\label{alg_1} 
\DontPrintSemicolon
\footnotesize
\KwIn{$Q \succeq 0$, $R\succ0$, $\gamma \in [0,1]$, $\ell \leftarrow 1$, $i^{e\ast(1)}$, $i^{eD}$, $\epsilon >0$, $\tau > 0$, $\underbar{$i$}_k^{e*}$, and $\beta_k$ for $k \in \mathbb N_{\rm E}$}
\While{$\ell \geq 2$ \text{and} $\|i^{e\ast(\ell+1)} - i^{e\ast(\ell)}\| \geq \tau$}
  {Construct $A(i^{e\ast(\ell)}), B(i^{e\ast(\ell)})$, and $\hat{x}_0$ in \eqref{h2pfm_ss_original}. Compute $K(i^{e\ast(\ell)})$ by \eqref{Klin_val},  $A_{cl}(i^{e\ast(\ell)})$ by \eqref{defAcl}, and $G_{i^{e\ast(\ell)}}(s)$ by \eqref{transfer_function}\\
  \For{$k \in \mathbb N_{\rm E}$}
 {
Compute $\Delta G_k$ in \eqref{defG_delta}. Let its system matrices be $\{A_{\Delta, k}, B_{\Delta, k}, C_{\Delta, k}\}$.
}
Define $\mathcal A$, $\mathcal B_1$, $\mathcal B_2$, and $\mathcal C$ in \eqref{defAAA}-\eqref{defCCC} \\
Compute the gradient of $J_{R}$ using (\ref{grad_obj})-\eqref{lyp}\\
Find optimal $\alpha_s$ using line search\\
Compute $i^{e\ast(\ell+1)}$ by \eqref{updateeq1}\\
Check if $i^{e\ast(\ell+1)}_k$ satisfy \eqref{ie_star} and update for $k \in \mathbb N_{\rm E}$, otherwise set to $\underbar{$i$}_k^{e\ast}$\\
Let $\ell\leftarrow \ell+1$\\
}
\KwOut{converged $i^{e\ast(\ell)}, K$}
\end{algorithm}

\subsection{Problem 2- Generalized EVCS integrated grid}\label{section_iii_C}
In our problem formulation so far, we have assumed that EVs only {\it consume} energy from the grid. However, at times of the day when the non-EV loads in the grid are high, the voltage stability indices at the buses, as defined in \eqref{vsi_main_exp}, may be low. For example, connecting multiple EVs working in the G2V mode to the end of a distribution line adversely affects its voltage profile, as shown in Fig. \ref{eig_time_vsi}(c). In such scenarios, EVs can also be used to {\it supply} energy to the grid in exchange of a fee, and help in recovering voltage stability. In other words, integrating the EVCSs with bidirectional capacity operating in vehicle-to-grid (V2G) mode can be an effective solution to maintain desired voltage profiles in the grid. In this section, we assess the impacts of unidirectional versus bidirectional EVCSs on charging control and incentivization, integrated with voltage stability analysis using VSIs. 

Let $\mathbb N_{\rm E_{uni}}$ and $\mathbb N_{\rm E_{bi}}$ denote the set of buses corresponding to uni- and bi-directional EVCSs, i.e., $\mathbb N_{\rm E_{uni}} \cup \mathbb N_{\rm E_{bi}} = \mathbb N_{\rm E}$ and $\mathbb N_{\rm E_{uni}} \cap \mathbb N_{\rm E_{bi}} = \emptyset$.
Similar to in section \ref{sec_proposed_alg}-A, we consider that the EV owners demand $i^{eD}$ as their charging rate, but the CPO offers $i^{e*}$ as the charging rate instead to improve the $\mathcal{H}_2$-norm of the grid transfer function, in exchange for incentives. In addition, the CPO now also wants to improve the VSI by buying power from the EVs using their bidirectional chargers, and, therefore, has to pay a certain fee to the corresponding EV owners. If the total spending budget of the CPO is fixed, this would lead to a compromise between improving the $\mathcal{H}_2$-norm and the VSI. We formulate this as an optimization problem by modifying the problem statement in \eqref{main_problem}-\eqref{ie_star} as follows:

\begin{align}
\min_{K,i^{e*}} J'= & (1-\gamma_1-\gamma_2) \underbrace{\int_{0}^{\infty} \|\hat{z}\|_2^2 dt}_{J_1} +\gamma_{1} \underbrace{\|i^{e*}-i^{eD}\|^2_{{2,\beta^{ip}}}}_{J_2} \nonumber \\
&+\gamma_{2}\underbrace{\|{\bf 1}_n-\mathcal{V}(i^{e*})\|_{\infty,\beta^{si}}}_{J_3}\label{main_problem2} \\
 \text{s.t} \quad & (\ref{overall_ss}), (\ref{LQR}) \; \mbox{hold},\\
& i_k^{e*} \in [\max\{\underbar{$i$}_k^{e*},0\}, {i}_k^{eD}], \quad k \in \mathbb N_{\rm E_{uni}},\label{ie_star_d1} \\ 
& i_k^{e*} \in [\underbar{$i$}_k^{e*}, {i}_k^{eD}], \quad k \in \mathbb N_{\rm E_{bi}},\label{ie_star_d2} 
\end{align}
where, $\gamma_1, \gamma_2 \in [0,\,1]$ are weighting constants that decide the priorities among the system $\mathcal{H}_2$-norm, VSI, and incentivization, $\beta^{ip}_k$ and $\beta^{si}_k$ are the price coefficients for $k \in \mathbb N_{\rm E}$, $\mathcal V := [\mathcal V_2, \ldots, \mathcal V_{n}]^{\top} \in \mathbb R^{n-1}$ with $\mathcal V_{h}$ as defined in \eqref{vsi_main_exp}. The VSI is included in the term $J_3$. 

\subsection{Solution Strategy of Problem 2}\label{section_iii_D}
Similarly to Section \ref{section_iii_B}, we consider 
\begin{align}
  \min_{i^{e*}} J'_{\rm R}= & (1-\gamma_1-\gamma_2)\underbrace{ \|G_{i^{e*}}\|_{\mathcal{H}_2}^2}_{J'_{R1}} + \gamma_1\underbrace{\|i^{e*}-i^{eD}\|^2_{2,\beta^{ip}}}_{J'_{R2}} \nonumber \\
   &+\gamma_2\underbrace{\|{\bf 1}_n - \mathcal{V}(i^{e*})\|_{\infty,\beta^{si}}}_{J'_{R3}} \label{relax_vsi} \\
\text{s.t.} \quad& \eqref{ie_star_d1}, \eqref{ie_star_d2}, \eqref{Klin}, \eqref{Klin_val}, \; \mbox{hold}, \label{relax_vsi_end}
\end{align}
where, $G_{i^{e*}}$ is defined in \eqref{transfer_function}. Similarly to \eqref{grad_obj}, the gradient of $J'_{R1}$, and $J'_{R2}$ are determined as 
\begin{equation}\label{grad_obj_2}
\scalebox{0.95}{$
    \left. \frac{\partial J'_{R1}}{\partial i^{e*}}\right|_{i^{e*}=i^{e*(\ell)}}= 2\mathcal B_2^{\top}L\mathcal B_1, ~
 \left. \frac{\partial J'_{R2}}{\partial i^{e\ast}}\right|_{i^{e*}=i^{e*(\ell)}}= 2\mathcal{D}^{ip} (i^{e\ast}- i^{eD}), $}
\end{equation}
where, $\mathcal{D}^{ip} := {\rm diag}(\beta^{ip}_k)_{k \in \mathbb N_{\rm E}}$. We next show how to obtain the gradient of $J'_{R3}$. 

Due to the complicated nonlinear relationship between $J'_{R3}(i^{e\ast})$ and $i^{e\ast}$, we consider a numerical approximation of the gradient. Considering \eqref{cl_state_1}, for $k \in \mathbb N_{\rm E}$ the $p_k$-th element of the gradient is expressed as 
\begin{equation}\label{JR3_grad}
   \left[\left.\frac{\partial J'_{R3}}{\partial i^{e}}\right|_{i^{e*}=i^{e*(\ell)}}\right]_k \approx \frac{J'_{R3}(i^{e*(\ell)}+e_k^p\epsilon) - J'_{R3}(i^{e*(\ell)})}{\epsilon}. 
\end{equation}
By repeating this procedure for every $k \in \mathbb N_{\rm E}$, we can obtain a numerical approximate of $\partial J'_{R3} / \partial i^{e*}$. The overall gradient of $J'_{R}$ is thereby approximated as
\begin{align} \label{J_prime_grad}
    &\hat{\nabla} J'_{R} := (1-\gamma_1-\gamma_2)2\mathcal B_2^{\top}L\mathcal B_1+\gamma_12\mathcal{D}^{ip} (i^{e\ast(\ell)}-i^{eD}) \nonumber \\
    &+\gamma_2 \left[\frac{J'_{R3}(i^{e*(\ell)}+e_k^p\epsilon) - J'_{R3}(i^{e*(\ell)})}{\epsilon}\right]_{k \in \mathbb N_{\rm E}}^{\top}
\end{align}
The optimization variable $i^{e*}$ is then updated as 
\begin{equation}\label{updateeq2}
    i^{e\ast(\ell+1)}=i^{e\ast(\ell)} - \alpha_v \hat{\nabla} J'_{R},
\end{equation}
where, $\alpha_v$ is the learning rate of the gradient. To find the optimal value of $\alpha_v$ the traditional Armijo-rule-based line search can be utilized. The entire procedure is summarized in Algorithm \ref{alg_2}.

\begin{algorithm}[h]
\caption{Algorithm for joint optimization and optimal control design \eqref{relax_vsi}-\eqref{relax_vsi_end}}\label{alg_2}
\DontPrintSemicolon
\footnotesize
\KwIn{$Q \succeq 0$, $R\succ0$, $\gamma_1, \gamma_2 \in [0,1]$, $\ell \leftarrow 1$, $i^{e*(1)}$ $ i^{eD}$, $\epsilon>0$, $\tau>0$, $\underbar{$i$}_k^{e*}$, $\beta_k^{ip}$, $\beta_k^{si}$ for $k \in \mathbb N_{\rm E}$}
\While{$\ell \geq 2$ \text{and} $\|i^{e\ast(\ell+1)} - i^{e\ast(\ell)}\| \geq \tau$}
{Construct $A(i^{e\ast(\ell)}), B(i^{e\ast(\ell)})$, and $\hat{x}_0$ in \eqref{h2pfm_ss_original}. Compute $K(i^{e\ast(\ell)})$ by \eqref{Klin_val},  $A_{cl}(i^{e\ast(\ell)})$ by \eqref{defAcl}, and $G_{i^{e\ast(\ell)}}(s)$ by \eqref{transfer_function}\\
  \For{$k \in \mathbb N_{\rm E}$}
 {
Compute $\Delta G_k$ in \eqref{defG_delta}. Let its system matrices be $\{A_{\Delta, k}, B_{\Delta, k}, C_{\Delta, k}\}$.
}
Define $\mathcal A$, $\mathcal B_1$, $\mathcal B_2$, and $\mathcal C$ in \eqref{defAAA}-\eqref{defCCC} \\
Compute the gradients of $ J'_{R1}$ and $J'_{R2}$ using \eqref{grad_obj_2}  \\
Import network load and line data and define EVCS bus index $\mathbb N_{\rm E_{\rm uni}}$, and $\mathbb N_{\rm E_{\rm bi}}$\\
Compute $\mathcal{V}_k$ for $k \in \{2,\ldots,n\}$ using \eqref{vsi_main_exp}, and $\mathcal V := [\mathcal V_2, \ldots, \mathcal V_n]^{\top}$. \\
Compute the gradient of $J'_{R3}$ using \eqref{JR3_grad}\\
Find optimal $\alpha_v$ using line search\\
Compute overall gradient of $J'_{R}$ using \eqref{J_prime_grad}\\
Update $i^{e*(\ell+1)}$ using \eqref{updateeq2}\\
Check if $i_k^{e*(\ell+1)}$ satisfy constraints \eqref{ie_star_d1}-\eqref{ie_star_d2} and update for $k\in \mathbb N_{\rm E}$, otherwise set to $\underbar{$i$}_k^{e*}$  \\

Let $\ell\leftarrow \ell+1$
}
\KwOut{converged $i^{e*(\ell)}, K$} 
\end{algorithm}

\section{Practical Implementation}\label{sec_prac_imp}
We next state the set of actions that the CPO and the EV owners need to take once the CPO runs Algorithm \ref{alg_1} or \ref{alg_2} in practice. Realistically speaking, every EV owner would want to submit their bid only a few hours before charging their EVs so that their bids are consistent with the latest pricing strategy for that day. At the same time, the CPO may need a sufficient amount of time to gather the bidding information from all interested customers so that Algorithms 1 or 2 can be run accurately. To strike a balance between these two factors, one idea can be to ask the customers to submit charging information one day ahead of time. This will include their charging energy demands and their desired power rate demands, which in most fast-charging or discharging EVCSs currently come in the form of discrete levels such as 50 kW, 75 kW, and 150 kW, as listed in \cite{eig_paper_5}. 

Using this information, the CPO constructs the vectors $E^{D}:=[E^D_{k}]_{k \in \mathbb N_{\rm E}}^{\top} \in\mathbb{R}^p$, and $P^{eD}:=[P_{k}^{eD}]_{k \in \mathbb N_{\rm E}}^{\top}\in\mathbb{R}^p$, and computes the charging or discharging current demand vectors using $i^{eD}:=P^{eD} \oslash v^{dc\ast} \in \mathbb R^p$, where, $v^{dc*}:= [v_k^{dc*}]_{k \in \mathbb N_{\rm E}}^{\top}$ and $v_k^{dc\ast}$ is the setpoint of $v_k^{dc}$ defined in \eqref{dc_bus}. For V2G mode, the sign of the charging power rate or charging current rate is considered to be negative, which is the opposite of the G2V mode. Using the chosen pricing co-efficient $\beta$ (in Problem 1) or $\beta^{si}$, and $\beta^{ip}$ (in problem 2), the CPO runs Algorithm \ref{alg_1} or \ref{alg_2} to determine the optimal charging current $i^{e*}$, and therefore, the optimal charging power rates $P^{e*}:= i^{e*} \odot v^{dc\ast} \in \mathbb R^p$. Additionally, the CPO will estimate the charging duration (in minutes) at the demanded charging rate or contracted discharging rate as 
$\mathcal{C}^{tD}:= 60 E^{D} \oslash P^{eD} \in\mathbb{R}^p$. Similarly, the charging or discharging duration at the optimal charging or discharging rate vector is calculated as $\mathcal{C}^{t*}:= 60E^{D} \oslash P^{e*} \in\mathbb{R}^p$. The charging price at the demanded charging or contracted discharging rate is calculated as $\mathcal{C}^{\mathcal{P}D}:=E^{D}\odot \beta$.
Based on this, the waiting time (in minutes), the dollar incentive, and the final charging cost are respectively calculated as $\mathcal{W}^{t}:=\mathcal{C}^{t*}-\mathcal{C}^{tD}$, $\mathcal{I}^{e}:=(\mathcal{C}^{\mathcal{P}D}\odot\mathcal{W}^{t})\oslash\mathcal{C}^{tD}$ and $\mathcal{C}^{\mathcal{P*}}:=\mathcal{C}^{\mathcal{P}D}-\mathcal{I}^{e}$, respectively.

The incentive offer is designed using a ``wait \& save'' strategy, similar to the protocols used for ride-sharing applications \cite{wait_save}. This information can be integrated into a mobile app, enabling direct communication between the CPO and EV owners. EV owners using the mobile app can access the recommended offers provided by the CPO and have the freedom to respond with their decisions. 

\begin{remark}
    Ideally, if any EV user rejects the incentivization offer, the CPO should re-optimize the remaining setpoints to still guarantee the minimum value of $J_R$ in \eqref{relax} or $J'_{R}$ in \eqref{relax_vsi}, and the new values of the waiting times and incentives should be communicated to the users who accepted their offers in the previous round. The user has the option of accepting or rejecting this new offer, based on which the CPO may need to re-optimize again. This negotiation can continue to the app until a consensus of acceptance or rejection is reached. For the sake of brevity, we skip this re-optimization scenario in this paper and reserve it for our future work.
\end{remark} 

\section {Numerical Results}\label{sec_results}
We validate our proposed method using the IEEE 33-bus radial feeder model. This model has an operating voltage of 12.66 kV, has 33 buses, and supports a 3.715 MW non-EV load. Traditional PI controller-based EVCSs are used with the dynamic model and power flow balance listed in Section II. We have seen in Fig. \ref{eig_time_vsi} that the damping worsens as the number of EVCSs increases. We aim to improve these damping factors using the co-design of $K$ and $i^{e\ast}$, as stated in Algorithm \ref{alg_1} or $K$ and $i^{e\ast}$ in Algorithm \ref{alg_2}. 
Please note that although the co-optimization is solved using the small-signal model, all simulations reported in the following subsections pertain to the nonlinear model of the 33-bus test system. All simulations are performed using MATLAB/Simulink.

\subsection{Unidirectional EVCSs}
\subsubsection{Single Port EVCSs}
Let us consider a charging scenario on any given day, where six EVs with the same battery capacity (75 kWh, 360V) submit the same charging demand $E^{D}=45 {\bf 1}_p$ kWh, i.e., 60\% of their capacity. We assume $p=3$ EVCSs case in Section \ref{sec_motiv}, i.e., $\mathbb N_{\rm E} = \{3,19,5\}$, each containing a single charging port. The customer-demanded charging power rates at those three EVCSs are $P^{eD}=[50,50,100]^{\top}$ kW, and corresponding charging current demands are $i^{eD}=[62.5,62.5,125]^{\top}$ A, considering $v^{dc*}=[800,800,800]^{\top}$ V. With this submitted information, the CPO runs Algorithm \ref{alg_1}, and divides the six EV loads equally among the three EVCSs, where each EVCS accommodates two EVs sequentially with no intervals of inactivity. 
To represent EV arrivals and departures in a condensed time frame under the state-of-charge (SOC) requirement, a scale-down of the demand needs to be done. We have scaled down the required 45 kWh demand to 0.04 kWh ({scale-down ratio of 1/1050}), and accordingly, the charging durations for the 50 kW and 100 kW chargers are reduced from their actual magnitudes of 54 minutes and 27 minutes to 2.8 seconds and 1.4 seconds of simulation time, respectively.

 \begin{figure}[h]
    \centering
     \includegraphics[scale=0.72]{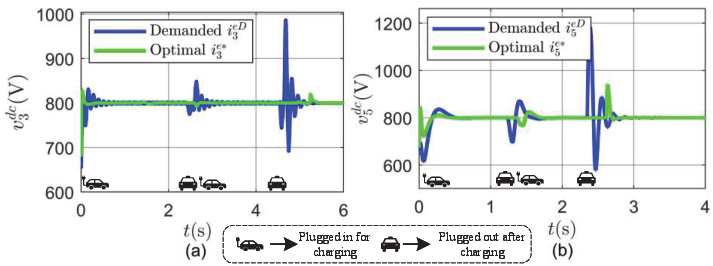}
  \caption{Comparison of DC bus voltage responses with (a) $i_{3}^{eD}=62.5$A, and $i_{3}^{e*}=55.97$A, and (b) $i_{5}^{eD}=125$A, and $i_{5}^{e*}=115.84$A.}
 \label{lqr_pi_h2_comp}
\end{figure}

With $\gamma=0$, Algorithm \ref{alg_1} yields the following optimal solutions: $i_{3}^{e\ast}=55.97$A, $i_{19}^{e\ast}=56.03$A, and $i_{5}^{e\ast}=115.84$A for the respective EVCSs with corresponding charging power rates  $P_{3}^{e\ast}=44.77$ kW, $P_{19}^{e\ast}=44.82$ kW, and $P_{5}^{e\ast}=92.67$ kW, respectively. The optimization reduces the $\mathcal{H}_2$-norm of the system transfer function $G_{i^{e\ast}}(s)$ in \eqref{transfer_function} from 0.0015 to $5.932 \times 10^{-4}$ with a computation time of 30.5 seconds. We apply the state-feedback controller $K$ and the optimal setpoints $i^{e\ast}$ in the original nonlinear model \eqref{overall_ss} of the 33-bus system, and consider a disturbance response due to the sudden connection of the EVs. Figure \ref{lqr_pi_h2_comp} shows the converter state responses (specifically, the DC bus voltage $v^{dc}$) for this charging scenario. Note that the voltage ripple in $v^{dc}$ depends on the current drawn. Therefore, the conventional PI-based model with the demanded $i^{eD}$ faces a higher risk of sustained voltage oscillations due to sudden plug-ins and outs of the high-power EV loads. This is also indicated in the voltage tolerance curves that are shown in \cite{EPRI}. 

A single figure is shown for the first two 50 kW charging stations as their capacities and optimization results are quite similar. Figure  \ref{lqr_pi_h2_comp} clearly illustrates that reducing the $\mathcal{H}_2$-norm leads to a significant improvement in transient response, as compared to the conventional PI-based controller with demand $i^{eD}$. The algorithm initially suggests a charging rate that is much lower than the customer's demand, potentially discouraging the customer's choice to use the EVCS. To address this, we introduce the incentivization term $J_{R2}$ in (\ref{relax}) that can be utilized by the CPO to ensure customer participation via a proper incentive. The CPO can tune the weighting factor $\gamma$ to attract EV customers. $\gamma=0$ prioritizes system dynamics over customer satisfaction, reducing the $\mathcal{H}_2$-norm from 0.0015 to $5.932 \times 10^{-4}$, while suggesting sub-demand charging currents. As $\gamma$ increases, both parts of the objective functions gain priority, leading to a trade-off.  As $\gamma$ approaches 1, emphasizing customer satisfaction over grid health, the algorithm suggests the use of charging currents that are closer to the customer demand but at the expense of degraded transient response of the grid states.

\subsubsection{Multi-port EVCSs}
\begin{figure}[h]
    \centering
    \includegraphics[scale=0.45]{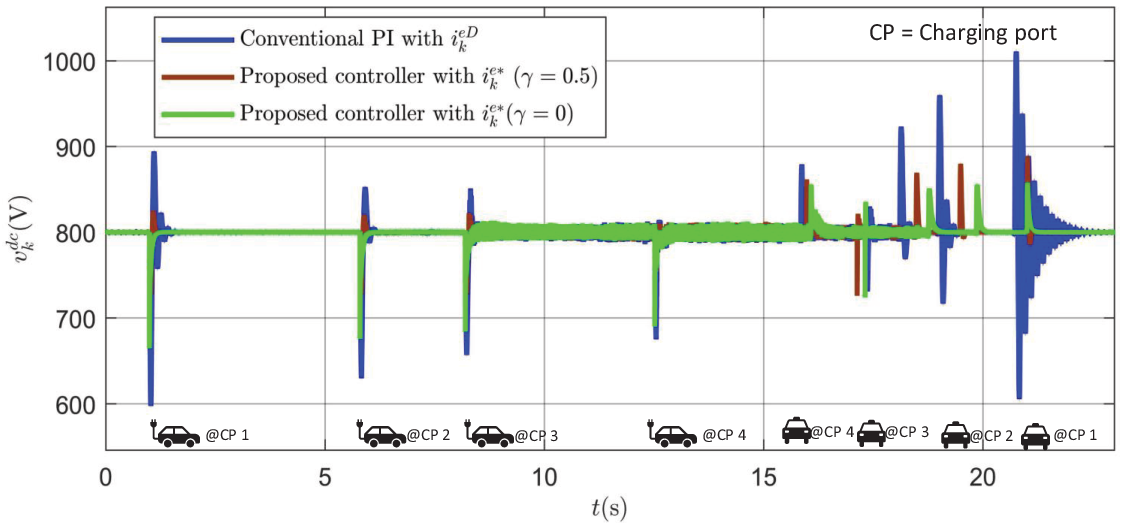}
    \caption{Comparison of DC bus voltages considering muti-port EVCSs. Here, the multi-port EVCS contains 4 charging ports (CPs), where EVs are plugged in and out at different times.}
    \label{h2_pi_multi}
\end{figure}

Moving beyond single charging port, we next study the impact of multi-port EVCSs on the closed-loop response. We consider the same three EVCSs in the 33-bus radial network, but now the first two EVCSs contain four charging ports each, and the third EVCS contains a single charging port. The charging rates for the first two EVCSs are 50 kW each, and that for the third EVCS is 100 kW. We assume the CPO allocates the first four EVs to the first EVCS and the remaining two to the third EVCS. Considering the condensed time frame, EV 1 joins the first charging port at the first EVCS at $t=1$s, and leaves at $t=20.7$s after receiving its 45 kWh of charging demand (i.e. 60\% of total SOC level of 75 kWh). The other EVs join and leave their respective charging ports while EV 1 is charging at the charging port 1 of the first EVCS. The transient response for the single charging port of the third EVCS is found to be similar to that in Fig. \ref{lqr_pi_h2_comp} (b), and therefore, is not repeated here. The responses for the multiport EVCS are shown in Fig. \ref{h2_pi_multi} using Algorithm \ref{alg_1}, for $\gamma=0$ and $\gamma=0.5$ cases. As expected, the transient response is best for $\gamma=0$ while the charging current is lower than the demand. The transient response worsens with increasing $\gamma$ while the CPO saves money by bringing the charging current closer to the demand.
\begin{table*}[ht]
\caption{Data table for negotiation between CPO and EV owner while utilizing the proposed incentivization strategy}
\centering
\scriptsize
\begin{tabular}{|c|c||cccc||cccc||ccc|}
\hline
\multirow{3}{*}{Cases} & \multirow{3}{*}{\begin{tabular}[c]{@{}c@{}}EVCS\\ ID\end{tabular}} & \multicolumn{4}{c|}{CPO Side}                                & \multicolumn{4}{c|}{Customer Side (Mobile App)}                             & \multicolumn{3}{c|}{Grid VSI Monitoring}  \\ \cline{3-13} 
                       &                                                                    & \multicolumn{1}{c|}{\multirow{2}{*}{$E^{D}$}} & \multicolumn{1}{c|}{\multirow{2}{*}{$P^{eD}$}} & \multicolumn{1}{c|}{\multirow{2}{*}{$P^{e*}$}} & \multirow{2}{*}{Decision} & \multicolumn{1}{c|}{\multirow{2}{*}{$\mathcal{C}^{t*},\mathcal{C}^{tD}$}} & \multicolumn{1}{c|}{\multirow{2}{*}{$\mathcal{W}^{t}$}} & \multicolumn{1}{c|}{\multirow{2}{*}{$\mathcal{I}^{e}$}} & \multirow{2}{*}{$\mathcal{C}^{\mathcal{P*}},\mathcal{C}^{\mathcal{P}D}$} & \multicolumn{1}{c|}{\multirow{2}{*}{$\mathcal{V}^{D}$}} & \multicolumn{1}{c|}{\multirow{2}{*}{$\mathcal{V}^{S}$}} & \multirow{2}{*}{$\mathcal{V}^{B}$} \\
                       &                                                                    & \multicolumn{1}{c|}{}                                   & \multicolumn{1}{c|}{}                                  & \multicolumn{1}{c|}{}                                  &                                        & \multicolumn{1}{c|}{}                                            & \multicolumn{1}{c|}{}                                            & \multicolumn{1}{c|}{}                                           &                                                      & \multicolumn{1}{c|}{}                                       & \multicolumn{1}{c|}{}                                       &                                      \\ \hline
\multirow{6}{*}{1}    & \multirow{2}{*}{EVCS-1}                                            & \multicolumn{1}{c|}{\multirow{2}{*}{45}}                & \multicolumn{1}{c|}{\multirow{2}{*}{50}}               & \multicolumn{1}{c|}{\multirow{2}{*}{45}}             & $\surd$                                & \multicolumn{1}{c|}{60}                                        & \multicolumn{1}{c|}{6}                                         & \multicolumn{1}{c|}{2}                                        & 16                                                 & \multicolumn{1}{c|}{-}                                       & \multicolumn{1}{c|}{0.9324}                                 & \multirow{2}{*}{0.936}               \\ \cline{6-12}
                       &                                                                    & \multicolumn{1}{c|}{}                                   & \multicolumn{1}{c|}{}                                  & \multicolumn{1}{c|}{}                                  & $\times$                               & \multicolumn{1}{c|}{54}                                          & \multicolumn{1}{c|}{0}                                           & \multicolumn{1}{c|}{0}                                          & 18                                                   & \multicolumn{1}{c|}{0.9316}                                 & \multicolumn{1}{c|}{-}                                       &                                      \\ \cline{2-13} 
                       & \multirow{2}{*}{EVCS-3}                                            & \multicolumn{1}{c|}{\multirow{2}{*}{45}}                & \multicolumn{1}{c|}{\multirow{2}{*}{100}}              & \multicolumn{1}{c|}{\multirow{2}{*}{90.8}}            & $\surd$                                & \multicolumn{1}{c|}{29.7}                                        & \multicolumn{1}{c|}{2.7}                                         & \multicolumn{1}{c|}{2.25}                                       & 20.25                                                & \multicolumn{1}{c|}{-}                                       & \multicolumn{1}{c|}{0.874}                                 & \multirow{2}{*}{0.884}               \\ \cline{6-12}
                       &                                                                    & \multicolumn{1}{c|}{}                                   & \multicolumn{1}{c|}{}                                  & \multicolumn{1}{c|}{}                                  & $\times$                               & \multicolumn{1}{c|}{27}                                          & \multicolumn{1}{c|}{0}                                           & \multicolumn{1}{c|}{0}                                          & 22.5                                                 & \multicolumn{1}{c|}{0.8729}                                 & \multicolumn{1}{c|}{-}                                       &                                      \\ \cline{2-13} 
                       & \multirow{2}{*}{EVCS-5}                                            & \multicolumn{1}{c|}{\multirow{2}{*}{45}}                & \multicolumn{1}{c|}{\multirow{2}{*}{150}}              & \multicolumn{1}{c|}{\multirow{2}{*}{139}}             & $\surd$                                & \multicolumn{1}{c|}{19.4}                                          & \multicolumn{1}{c|}{1.4}                                           & \multicolumn{1}{c|}{1.75}                                        & 20.75                                                  & \multicolumn{1}{c|}{-}                                       & \multicolumn{1}{c|}{0.8001}                                 & \multirow{2}{*}{0.810}               \\ \cline{6-12}
                       &                                                                    & \multicolumn{1}{c|}{}                                   & \multicolumn{1}{c|}{}                                  & \multicolumn{1}{c|}{}                                  & $\times$                               & \multicolumn{1}{c|}{18}                                          & \multicolumn{1}{c|}{0}                                           & \multicolumn{1}{c|}{0}                                          & 22.5                                                 & \multicolumn{1}{c|}{0.7988}                                 & \multicolumn{1}{c|}{-}                                       &                                      \\ \hline
\multirow{6}{*}{2}    & \multirow{2}{*}{EVCS-1}                                            & \multicolumn{1}{c|}{\multirow{2}{*}{45}}                & \multicolumn{1}{c|}{\multirow{2}{*}{30}}               & \multicolumn{1}{c|}{\multirow{2}{*}{25.5}}               & $\surd$                                & \multicolumn{1}{c|}{106}                                         & \multicolumn{1}{c|}{16}                                          & \multicolumn{1}{c|}{4}                                        & 18.5                                                   & \multicolumn{1}{c|}{-}                                       & \multicolumn{1}{c|}{0.9341}                                 & \multirow{2}{*}{0.936}               \\ \cline{6-12}
                       &                                                                    & \multicolumn{1}{c|}{}                                   & \multicolumn{1}{c|}{}                                  & \multicolumn{1}{c|}{}                                  & $\times$                               & \multicolumn{1}{c|}{90}                                          & \multicolumn{1}{c|}{0}                                           & \multicolumn{1}{c|}{0}                                          & 22.5                                                 & \multicolumn{1}{c|}{0.9328}                                 & \multicolumn{1}{c|}{-}                                       &                                      \\ \cline{2-13} 
                       & \multirow{2}{*}{EVCS-3}                                            & \multicolumn{1}{c|}{\multirow{2}{*}{45}}                & \multicolumn{1}{c|}{\multirow{2}{*}{75}}              & \multicolumn{1}{c|}{\multirow{2}{*}{67.1}}            & $\surd$                                 & \multicolumn{1}{c|}{40.2}                                          & \multicolumn{1}{c|}{4.2}                                         & \multicolumn{1}{c|}{3.2}                                          & 23.8                                                 & \multicolumn{1}{c|}{-}                                       & \multicolumn{1}{c|}{0.8768}                                 & \multirow{2}{*}{0.884}               \\ \cline{6-12}
                       &                                                                    & \multicolumn{1}{c|}{}                                   & \multicolumn{1}{c|}{}                                  & \multicolumn{1}{c|}{}                                  & $\times$                               & \multicolumn{1}{c|}{36}                                        & \multicolumn{1}{c|}{0}                                           & \multicolumn{1}{c|}{0}                                          & 23.8                                                   & \multicolumn{1}{c|}{0.8754}                                 & \multicolumn{1}{c|}{-}                                       &                                      \\ \cline{2-13} 
                       & \multirow{2}{*}{EVCS-5}                                            & \multicolumn{1}{c|}{\multirow{2}{*}{45}}                & \multicolumn{1}{c|}{\multirow{2}{*}{120}}               & \multicolumn{1}{c|}{\multirow{2}{*}{107.4}}             & $\surd$                                & \multicolumn{1}{c|}{25.1}                                          & \multicolumn{1}{c|}{2.6}                                           & \multicolumn{1}{c|}{3.1}                                          & 23.9                                                   & \multicolumn{1}{c|}{-}                                       & \multicolumn{1}{c|}{0.8027}                                 & \multirow{2}{*}{0.810}               \\ \cline{6-12}
                       &                                                                    & \multicolumn{1}{c|}{}                                   & \multicolumn{1}{c|}{}                                  & \multicolumn{1}{c|}{}                                  & $\times$                                & \multicolumn{1}{c|}{22.5}                                          & \multicolumn{1}{c|}{0}                                           & \multicolumn{1}{c|}{0}                                          & 27                                                   & \multicolumn{1}{c|}{0.8015}                                 & \multicolumn{1}{l|}{-}                                       &                                      \\ \hline
\end{tabular}\label{datalog_table}\vspace{1mm}\\
\footnotesize{\footnotesize EV owner decision: $\surd=$Accept, $\times=$Reject}. Unit: $\mathcal{C}^{tD}, \mathcal{C}^{t*},\mathcal{W}^{t}$ in minutes, $E^{D}$ in kWh, $P^{eD}$, $P^{e\ast}$ in kW, and $\mathcal{I}^{e}$, $\mathcal{C}^{\mathcal{P}D}$, $\mathcal{C}^{\mathcal{P\ast}}$ in \$. \\
\end{table*}
\subsubsection{Negotiation through Incentivization}
 
For each $k \in \mathbb{N}_{\text{E}}$, we set $\beta_k = 0.4$/kWh during off-peak and $\beta_k = 0.5$/kWh during peak hours (for charging power rates under 50 kW), and $\beta_k = 0.5$/kWh during off-peak and $\beta_k = 0.6$/kWh during peak hours (for charging rates exceeding 50 kW), based on average prices from US charging providers \cite{charging_price}.
We assume $p=3$ EVCSs connected to buses 3, 5, and 15, i.e., $\mathbb{N}_{\rm E}:=\{3,5,15\}$, each with a charging energy demand of 45 kWh, as shown in Table \ref{datalog_table}. Here, the EVCS ID is the corresponding serial number of the EVCSs as shown in Fig.~\ref{33_bus}. The optimal charging power rates $P^{e*}$ from Algorithm 1 are lower than $P^{eD}$, extending the charging duration $\mathcal{C}^{tD}$ (Section IV). We introduce the ``wait \& save'' strategy \cite{wait_save}, where EV customers trade a brief waiting time $\mathcal{W}^{t}$ for saving money, reducing their final charging price $\mathcal{C}^{\mathcal{P*}}$. To explain this, we consider two cases: Case-1, EV charging loads in the mentioned EVCSs with three different demands $P^{eD}$ in off-peak, and Case-2, EV charging loads with lower demands $P^{eD}$ in peak hours.

\begin{figure}[h]
    \centering
    \includegraphics[width=0.5\textwidth]{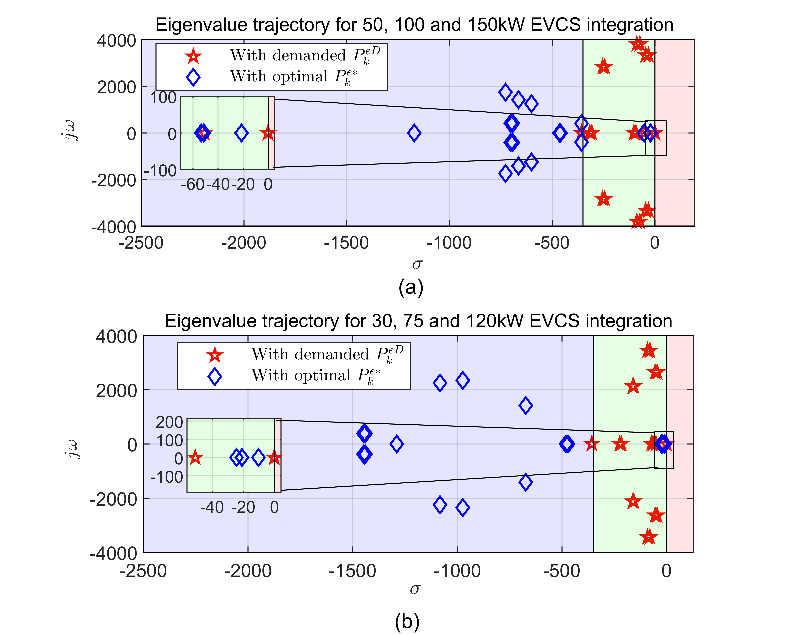}
    \caption{Eigenvalue plots showing improvement of closed-loop damping factors for the dominant poles for cases 1 and 2 of Table~\ref{datalog_table}. }
    \label{eigen_traj_all}
\end{figure}

\begin{figure}[h]
    \centering
    \includegraphics[scale=0.61]{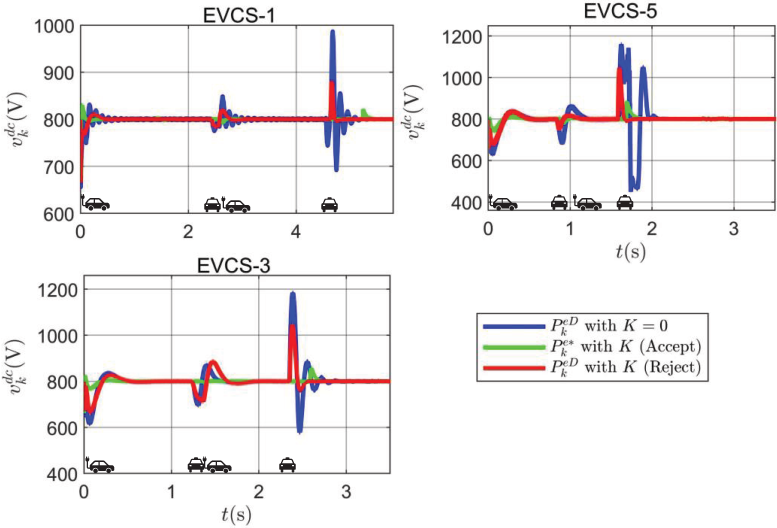}
    \caption{DC bus voltages of EVCSs 1, 3, and 5, for Case 1 of Table~\ref{datalog_table}.}
    \label{dc_bus_demo}
\end{figure}

The corresponding data log observed from the CPO side and on the mobile app screen of the customer side are shown in Table \ref{datalog_table}. Additionally, the VSI at the EVCS buses, i.e. $k\in\mathbb{N}_{\rm E}$ is computed following \cite{Chakra}, and tabulated in Table \ref{datalog_table}, where, the base load VSI (i.e. with no EV load at EVCS buses) is denoted as $\mathcal{V}^{\mathcal{B}}$ at $k\in\mathbb{N}_{\rm L}$. In case 1, three EV loads join at the three EVCS buses with $P^{eD}$ of 50 kW, 100 kW, and 150 kW, respectively. Thereafter, the CPO runs Algorithm 1, and determines the optimal $P^{e*}$ as 45 kW, 90.8 kW, and 139 kW, respectively. Therefore, the charging durations $\mathcal{C}^{t*}$ become 60 minutes, 29.7 minutes, and 19.4 minutes, respectively. The charging durations $\mathcal{C}^{tD}$ if the EVs stick to their demanded power, however, are 54 minutes, 27 minutes, and 18 minutes, respectively. The respective customers, therefore, will have to wait for an extra 6 minutes, 2.7 minutes, and 1.4 minutes, which eventually saves them \$2, \$2.25, and \$1.75 as incentives $\mathcal{I}^{e}$ (calculated based on the equivalent dollar value for staying the extra time), helping to reduce their final charging price $\mathcal{C}^{\mathcal{P*}}$. The utilization of lower $P^{e*}$ also slightly improves the VSI than $P^{eD}$, as shown in Table \ref{datalog_table}, where $\mathcal{V}^{\mathcal{S}}$ is the VSI when the customer agrees to go with optimal charging rate, and $\mathcal{V}^{D}$ is the VSI when they reject the idea of the optimal rate and prefer to stick to their demand. The same observation is found for case 2, as shown in Table \ref{datalog_table}. 
The closed-loop eigenvalue plots for cases 1 \& 2, and the respective transient responses of the DC bus voltage are shown in Figs. \ref{eigen_traj_all} and \ref{dc_bus_demo}, respectively. The dynamic response of the voltage is significantly better if the consumer accepts the proposed charging rate. If the consumer rejects the suggested charging rate, the closed-loop response degrades but still remains better than that for the demanded charging rate. In summary, the proposed dynamics-aware charging scheme provides beneficial returns for both the grid operator and the EV customers. Concerns on voltage stability can be overcome by utilizing both the G2V and V2G modes. This is discussed in the next subsection.
\subsection{Both Unidirectional and Bidirectional EVCSs}
For validating Problem 2, we consider the same IEEE 33 bus model as shown in Fig.~\ref{33_bus} with $p=10$ EVCSs, where 7 EVCSs connected at buses 3, 5, 9, 19, 21, 26, and 32 are considered to be operating in uni-directional G2V mode, and 3 EVCSs at buses 11, 15, and 17 are assumed to operate as bi-directional EVCSs. The nominal VSI without any EVCS in the network (as shown in Fig. \ref{eig_time_vsi}(c)) is found to be the lowest at bus 18 with a value of 0.7141. This VSI is marked as $\mathcal{V}_{\rm min}$. Due to the radial nature of the network, it can be easily seen that the VSIs at all other buses will improve, in general, if  $\mathcal{V}_{\rm min}$ improves.
We consider that two EVs in each EVCS with the same battery capacity (120 kWh, 360V) submit their charging demand as 90 kWh, due to which the VSI of the critical bus reduces to $\mathcal{V}_{\rm min}=0.5768$. The CPO, therefore, decides to purchase power from the EVCSs that have bi-directional capabilities. The demanded charging and contracted discharging rates, denoted as $P^{eD}$, are shown in Table \ref{bidirectional_table}.
The bracketed sign indicates the EVCSs that operate in V2G mode only. Given this submitted demand information, the CPO runs a power flow program and obtains the critical VSI (without applying Algorithm 2) as $\mathcal{V}_{\rm min}=0.654$, indicating an improvement in voltage due to the added energy supplied by the EVCSs. However, the $\mathcal{H}_2$-norm of the system is still sub-optimal as the charging/discharging current setpoints are not optimized. The CPO, therefore, runs Algorithm 2 to co-optimize $\mathcal{H}_2$-norm, VSI, and the incentives. 

\begin{table}[h]
\scriptsize
\caption{Results for IEEE 33-bus test model with 10 EVCSs with bidirectional capability}
\begin{tabular}{|cccccccc|}
\hline
\multicolumn{1}{|c|}{ID}      & \multicolumn{1}{c|}{$P^{eD}$} & \multicolumn{1}{c|}{$P^{e*}$} & \multicolumn{1}{c|}{$\mathcal{W}^t$} & \multicolumn{1}{c|}{$\mathcal{I}^e$} & \multicolumn{1}{c|}{$\sum\mathcal{I}^e$}       & \multicolumn{1}{c|}{ $\|G_{i^{e*}}\|_{\mathcal{H}_2}^{2}$} & $\mathcal{V}_{\rm min}$          \\ \hline

\multicolumn{8}{|c|}{Case 1: $\gamma_1=0$, $\gamma_2=0$}               \\ \hline
\multicolumn{1}{|c|}{1}  & \multicolumn{1}{c|}{50}       & \multicolumn{1}{c|}{42.5}     & \multicolumn{1}{c|}{19.06}           & \multicolumn{1}{c|}{6.35}          & \multicolumn{1}{c|}{\multirow{10}{*}{51.47}} & \multicolumn{1}{c|}{\multirow{10}{*}{\parbox{1cm}{\hspace{0.1cm}9.65e-4\\\ \{9.71e-4\}}}}  & \multirow{10}{*}{0.630} \\ \cline{1-5}
\multicolumn{1}{|c|}{2}  & \multicolumn{1}{c|}{100}      & \multicolumn{1}{c|}{85}       & \multicolumn{1}{c|}{9.53}            & \multicolumn{1}{c|}{7.94}          & \multicolumn{1}{c|}{}                        & \multicolumn{1}{c|}{}                                  &                          \\ \cline{1-5}
\multicolumn{1}{|c|}{3}  & \multicolumn{1}{c|}{100}      & \multicolumn{1}{c|}{86.85}    & \multicolumn{1}{c|}{8.18}            & \multicolumn{1}{c|}{6.82}          & \multicolumn{1}{c|}{}                        & \multicolumn{1}{c|}{}                                  &                          \\ \cline{1-5}
\multicolumn{1}{|c|}{4}  & \multicolumn{1}{c|}{(175)}     & \multicolumn{1}{c|}{(164.88)}  & \multicolumn{1}{c|}{1.89}            & \multicolumn{1}{c|}{2.75}          & \multicolumn{1}{c|}{}                        & \multicolumn{1}{c|}{}                                  &                          \\ \cline{1-5}
\multicolumn{1}{|c|}{5}  & \multicolumn{1}{c|}{(150)}     & \multicolumn{1}{c|}{(143.15)}  & \multicolumn{1}{c|}{1.72}            & \multicolumn{1}{c|}{2.15}          & \multicolumn{1}{c|}{}                        & \multicolumn{1}{c|}{}                                  &                          \\ \cline{1-5}
\multicolumn{1}{|c|}{6}  & \multicolumn{1}{c|}{(175)}     & \multicolumn{1}{c|}{(169.83)}  & \multicolumn{1}{c|}{0.94}            & \multicolumn{1}{c|}{1.37}          & \multicolumn{1}{c|}{}                        & \multicolumn{1}{c|}{}                                  &                          \\ \cline{1-5}
\multicolumn{1}{|c|}{7}  & \multicolumn{1}{c|}{50}       & \multicolumn{1}{c|}{42.5}     & \multicolumn{1}{c|}{19.06}           & \multicolumn{1}{c|}{6.35}          & \multicolumn{1}{c|}{}                        & \multicolumn{1}{c|}{}                                  &                          \\ \cline{1-5}
\multicolumn{1}{|c|}{8}  & \multicolumn{1}{c|}{100}      & \multicolumn{1}{c|}{85.94}    & \multicolumn{1}{c|}{8.83}            & \multicolumn{1}{c|}{7.36}          & \multicolumn{1}{c|}{}                        & \multicolumn{1}{c|}{}                                  &                          \\ \cline{1-5}
\multicolumn{1}{|c|}{9}  & \multicolumn{1}{c|}{50}       & \multicolumn{1}{c|}{42.5}     & \multicolumn{1}{c|}{19.06}           & \multicolumn{1}{c|}{6.35}          & \multicolumn{1}{c|}{}                        & \multicolumn{1}{c|}{}                                  &                          \\ \cline{1-5}
\multicolumn{1}{|c|}{10} & \multicolumn{1}{c|}{100}      & \multicolumn{1}{c|}{91.77}    & \multicolumn{1}{c|}{4.84}            & \multicolumn{1}{c|}{4.03}          & \multicolumn{1}{c|}{}                        & \multicolumn{1}{c|}{}                                  &                          \\ \hline

\multicolumn{8}{|c|}{Case 2: $\gamma_1=0.33$, $\gamma_2=0.33$}          \\ \hline
\multicolumn{1}{|c|}{1}  & \multicolumn{1}{c|}{50}       & \multicolumn{1}{c|}{45.87}    & \multicolumn{1}{c|}{9.72}            & \multicolumn{1}{c|}{3.24}          & \multicolumn{1}{c|}{\multirow{10}{*}{19.37}} & \multicolumn{1}{c|}{\multirow{10}{*}{\parbox{1cm}{\hspace{0.2cm}1.5e-3\\\{1.53e-3\}}}}           & \multirow{10}{*}{0.691} \\ \cline{1-5}
\multicolumn{1}{|c|}{2}  & \multicolumn{1}{c|}{100}      & \multicolumn{1}{c|}{94.83}    & \multicolumn{1}{c|}{2.94}            & \multicolumn{1}{c|}{2.45}          & \multicolumn{1}{c|}{}                        & \multicolumn{1}{c|}{}                                  &                          \\ \cline{1-5}
\multicolumn{1}{|c|}{3}  & \multicolumn{1}{c|}{100}      & \multicolumn{1}{c|}{95.91}    & \multicolumn{1}{c|}{2.30}            & \multicolumn{1}{c|}{1.92}          & \multicolumn{1}{c|}{}                        & \multicolumn{1}{c|}{}                                  &                          \\ \cline{1-5}
\multicolumn{1}{|c|}{4}  & \multicolumn{1}{c|}{(175)}     & \multicolumn{1}{c|}{(171.83)}  & \multicolumn{1}{c|}{0.55}            & \multicolumn{1}{c|}{0.80}          & \multicolumn{1}{c|}{}                        & \multicolumn{1}{c|}{}                                  &                          \\ \cline{1-5}
\multicolumn{1}{|c|}{5}  & \multicolumn{1}{c|}{(150)}     & \multicolumn{1}{c|}{(147.84)}  & \multicolumn{1}{c|}{0.53}            & \multicolumn{1}{c|}{0.66}          & \multicolumn{1}{c|}{}                        & \multicolumn{1}{c|}{}                                  &                          \\ \cline{1-5}
\multicolumn{1}{|c|}{6}  & \multicolumn{1}{c|}{(175)}     & \multicolumn{1}{c|}{(173.38)}  & \multicolumn{1}{c|}{0.29}            & \multicolumn{1}{c|}{0.42}          & \multicolumn{1}{c|}{}                        & \multicolumn{1}{c|}{}                                  &                          \\ \cline{1-5}
\multicolumn{1}{|c|}{7}  & \multicolumn{1}{c|}{50}       & \multicolumn{1}{c|}{44.98}    & \multicolumn{1}{c|}{12.05}           & \multicolumn{1}{c|}{4.02}          & \multicolumn{1}{c|}{}                        & \multicolumn{1}{c|}{}                                  &                          \\ \cline{1-5}
\multicolumn{1}{|c|}{8}  & \multicolumn{1}{c|}{100}      & \multicolumn{1}{c|}{95.63}    & \multicolumn{1}{c|}{2.47}            & \multicolumn{1}{c|}{2.06}          & \multicolumn{1}{c|}{}                        & \multicolumn{1}{c|}{}                                  &                          \\ \cline{1-5}
\multicolumn{1}{|c|}{9}  & \multicolumn{1}{c|}{50}       & \multicolumn{1}{c|}{46.63}    & \multicolumn{1}{c|}{7.81}            & \multicolumn{1}{c|}{2.60}          & \multicolumn{1}{c|}{}                        & \multicolumn{1}{c|}{}                                  &                          \\ \cline{1-5}
\multicolumn{1}{|c|}{10} & \multicolumn{1}{c|}{100}      & \multicolumn{1}{c|}{97.42}    & \multicolumn{1}{c|}{1.43}             & \multicolumn{1}{c|}{1.20}          & \multicolumn{1}{c|}{}                        & \multicolumn{1}{c|}{}                                  &                          \\ \hline

\multicolumn{8}{|c|}{Case 3:  $\gamma_1=0$, $\gamma_2=0.6$}            \\ \hline
\multicolumn{1}{|c|}{1}  & \multicolumn{1}{c|}{50}       & \multicolumn{1}{c|}{42.75}    & \multicolumn{1}{c|}{18.32}           & \multicolumn{1}{c|}{6.10}          & \multicolumn{1}{c|}{\multirow{10}{*}{37.13}} & \multicolumn{1}{c|}{\multirow{10}{*}{\parbox{1cm}{\hspace{0.2cm}1.6e-3\\\{1.64e-3\}}}}          & \multirow{10}{*}{0.770} \\ \cline{1-5}
\multicolumn{1}{|c|}{2}  & \multicolumn{1}{c|}{100}      & \multicolumn{1}{c|}{92.11}    & \multicolumn{1}{c|}{4.63}            & \multicolumn{1}{c|}{3.86}          & \multicolumn{1}{c|}{}                        & \multicolumn{1}{c|}{}                                  &                          \\ \cline{1-5}
\multicolumn{1}{|c|}{3}  & \multicolumn{1}{c|}{100}      & \multicolumn{1}{c|}{92.73}    & \multicolumn{1}{c|}{4.23}            & \multicolumn{1}{c|}{3.53}          & \multicolumn{1}{c|}{}                        & \multicolumn{1}{c|}{}                                  &                          \\ \cline{1-5}
\multicolumn{1}{|c|}{4}  & \multicolumn{1}{c|}{(175)}     & \multicolumn{1}{c|}{(168.38)}  & \multicolumn{1}{c|}{1.21}            & \multicolumn{1}{c|}{1.76}          & \multicolumn{1}{c|}{}                        & \multicolumn{1}{c|}{}                                  &                          \\ \cline{1-5}
\multicolumn{1}{|c|}{5}  & \multicolumn{1}{c|}{(150)}     & \multicolumn{1}{c|}{(144.11)}  & \multicolumn{1}{c|}{1.47}            & \multicolumn{1}{c|}{1.84}          & \multicolumn{1}{c|}{}                        & \multicolumn{1}{c|}{}                                  &                          \\ \cline{1-5}
\multicolumn{1}{|c|}{6}  & \multicolumn{1}{c|}{(175)}     & \multicolumn{1}{c|}{(169.45)}  & \multicolumn{1}{c|}{1.01}            & \multicolumn{1}{c|}{1.47}          & \multicolumn{1}{c|}{}                        & \multicolumn{1}{c|}{}                                  &                          \\ \cline{1-5}
\multicolumn{1}{|c|}{7}  & \multicolumn{1}{c|}{50}       & \multicolumn{1}{c|}{42.5}     & \multicolumn{1}{c|}{19.06}           & \multicolumn{1}{c|}{6.35}          & \multicolumn{1}{c|}{}                        & \multicolumn{1}{c|}{}                                  &                          \\ \cline{1-5}
\multicolumn{1}{|c|}{8}  & \multicolumn{1}{c|}{100}      & \multicolumn{1}{c|}{92.53}    & \multicolumn{1}{c|}{4.36}            & \multicolumn{1}{c|}{3.63}          & \multicolumn{1}{c|}{}                        & \multicolumn{1}{c|}{}                                  &                          \\ \cline{1-5}
\multicolumn{1}{|c|}{9}  & \multicolumn{1}{c|}{50}       & \multicolumn{1}{c|}{43.24}    & \multicolumn{1}{c|}{16.88}           & \multicolumn{1}{c|}{5.63}          & \multicolumn{1}{c|}{}                        & \multicolumn{1}{c|}{}                                  &                          \\ \cline{1-5}
\multicolumn{1}{|c|}{10} & \multicolumn{1}{c|}{100}      & \multicolumn{1}{c|}{93.81}    & \multicolumn{1}{c|}{3.56}            & \multicolumn{1}{c|}{2.96}          & \multicolumn{1}{c|}{}                        & \multicolumn{1}{c|}{}                                  &                          \\ \hline
\end{tabular}\label{bidirectional_table}\\
\end{table}
The results from Algorithm 2 are presented in Table \ref{bidirectional_table}, where `ID' represents the identification number of EVCSs as marked in Fig. \ref{33_bus}. For Case 1, we assume $\gamma_1=\gamma_2=0$, giving higher priority to improve the system $\mathcal{H}_2$-norm. In this case, the algorithm provides lower charging and discharging power rates, denoted as $P^{e*}$ in the table for the corresponding EVCSs, ensuring a lower $\mathcal{H}_2$-norm of 9.65e-4. The computational time is 155.29 seconds. The optimization also ensures an improvement in VSI (i.e., 0.630) from the unidirectional non-optimal case (i.e., 0.5768). However, this costs the CPO an incentive of \$51.47. In case 3, the CPO aims to reduce the incentive and give more priority to the VSI improvement by setting $\gamma_1=0$ and $\gamma_2=0.6$. The algorithm provides a slightly higher power rate by sacrificing the $\mathcal{H}_2$-norm from $9.65e{-4}$ to $1.6e{-3}$, but ensures a lower incentive of \$37.13 and a higher VSI of 0.770. In case 2, equal priority is given to $\mathcal{H}_2$-norm improvement, customer satisfaction, and VSI improvement. In this case, the charging and discharging power rates are closer to customer demand, which reduces incentive costs, but both the system $\mathcal{H}_2$-norm and VSI suffer a bit. 

\begin{figure}[h]
    \centering
    \includegraphics[width=0.45\textwidth]{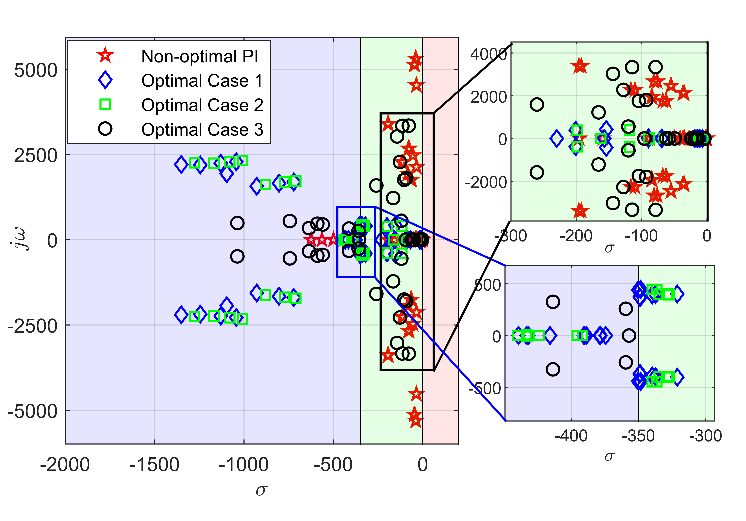}
    \caption{Eigenvalue comparison between four models: one with only PI control, and the ones with both PI control and LQR control for the three cases of Table \ref{bidirectional_table}. Case 1: $\gamma_1=0$, $\gamma_2=0$, case 2: $\gamma_1=0$, $\gamma_2=0.6$, and case 3: $\gamma_1=0.33$, $\gamma_2=0.33$.}
    \label{eig_bi}
\end{figure}

\begin{figure}[h]
    \centering
    \includegraphics[width=0.45\textwidth]{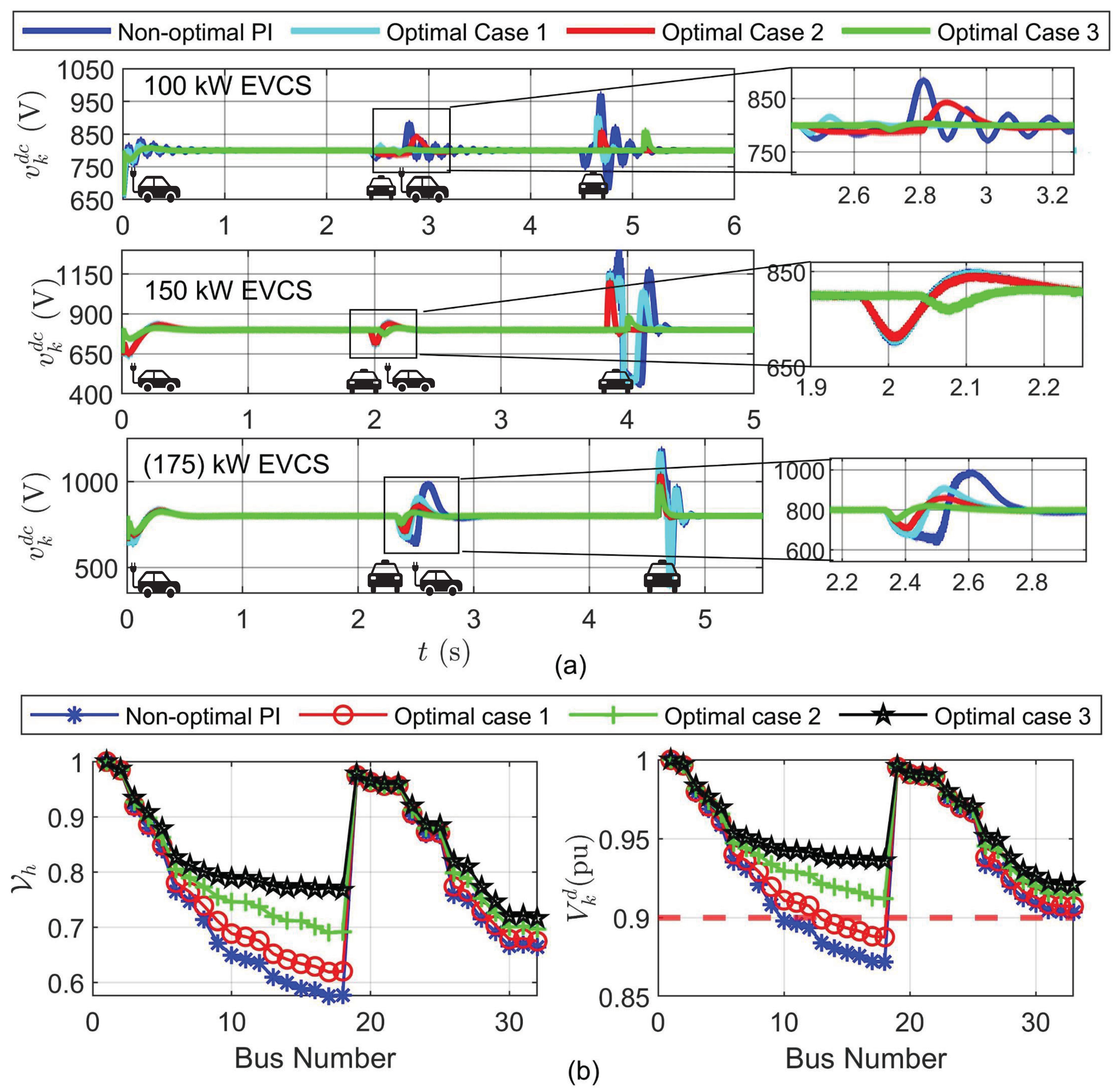}
    \caption{(a) DC bus voltage responses for all cases of Table \ref{bidirectional_table} versus the PI-only model, and (b) Trends of VSI values for PI-only control versus the PI+LQR control with different weighting factors.}
    \label{time_vsi}
\end{figure}

To test the robustness of our proposed method, we run Algorithm 2 using the nominal model of the test system but implement the solution considering that the system has undergone a 30\% increase in every non-EV load (i.e., $P_k^L, Q_k^L$). The resulting values of the $\mathcal{H}_2$-norm are listed within $\{.\}$ in column 7 of Table \ref{bidirectional_table}. It can be seen that the $\mathcal{H}_2$-norm values are only 0.62\%, 2\%, and 2.5\% higher than their respective nominal values, which verifies the low sensitivity of Algorithm 2 to model uncertainty.
The eigenvalue shift shown in Fig. \ref{eig_bi} also suggests an improvement in the small-signal dynamic performance. The time-domain responses of the DC bus voltages for selected EVCSs are shown in Fig. \ref{time_vsi}(a). For Case 1, where $\gamma_1=0=\gamma_2=0$, the transient response of the bus voltage is significantly better than the other cases. For Case 3, the system $\mathcal{H}_2$-norm is reduced at the expense of maximizing the VSI and reducing the incentive. The DC bus voltage responses show this degradation, with the corresponding closed-loop eigenvalues marked as green boxes in Fig. \ref{eig_bi} moving closer to the imaginary axis. A similar and consistent observation is made for case 2 as well. The trends of the VSI and the bus voltages using Algorithm 2 for the above three cases are shown in Fig. \ref{time_vsi}(b). The figure indicates how the CPO may decide on the interplay between the three objectives in (\ref{main_problem2}) to ensure a balance between grid health, EV customer satisfaction, and total expenditure. 

\section{Conclusions}\label{conclusion}
We presented a joint optimization and optimal control design for charging control of EVs in power distribution networks that benefits the grid operator by improving the small-signal damping performance of the grid states, and also the EV owners by providing them with financial incentives in exchange for slightly higher waiting times than what they intended. The damping improvement is achieved by a joint LQR state-feedback design for the charging current and by reducing the $\mathcal{H}_2$-norm of the system transfer function by finding the optimal setpoints for these charging currents. Both unidirectional and bidirectional charging scenarios are considered. The latter is seen to bring out interesting trade-offs between small-signal stability, voltage stability, and incentivization. The proposed approach is validated using simulations on a benchmark distribution grid, showcasing its effectiveness in mitigating voltage oscillations as well as in optimizing the charging costs of EVs. One challenge for implementing the proposed controller on a massive scale can be the associated peer-to-peer communication cost needed for the state feedback between the EVCSs. One potential solution for this issue can be to divide the network into zones and impose a block-diagonally dominant structure on $K$ so that the operator would need to install communication links only within the zones. We reserve such zonal analyses for our future work. Detailed cost-benefit analyses also need to be performed to verify the justification of the proposed incentivization method as an {\it insurance} against grid instability, and how that correlates with temporal scheduling of EVs and demand response. 

\section*{Appendix}
{\textit{A. EVCS Parameters} \\
 The magnitudes of the model parameters used for the design of the 50 kW EVCS are summarized as follows: For any $k \in \mathbb N_{\rm E}$, $L_k^g=L_k^c=2$ mH, $C_k^f=30$ $\mu$F, $v_k^{dc}=800$V, $\kappa_{k}^{P1}=1.71$, $\kappa_{k}^{I1}=672.66$, $\kappa_{k}^{P2}=0.5$, $\kappa_{k}^{I2}=5$,  $\kappa_{k}^{P3}=25$, $\kappa_{k}^{I3}=500$, $\kappa_{k}^{P4}=25$, $\kappa_{k}^{I4}=500$, $C^{dc}_k=5600$ $\mu$F.
 
\textit{B. Derivation for determining the gradient of $J_{R}$ in (\ref{relax})} \\
Following \cite{h2pfm}, the gradient of $J_{R1}$ part of $J_{R}$ can be expressed as:
\begin{equation}\nonumber
J_{R1}=\|G_{i^{e*}}\|_{\mathcal{H}_2}^2=\mathrm{tr}\left(\mathcal{B}^{\top}L{\mathcal{B}}\right)
\end{equation}
\begin{equation} \nonumber
=\mathrm{tr}\left(({\mathcal{B}_1}+{\mathcal{B}_2\epsilon {\bf 1}_p)}^{\top}L({\mathcal{B}_1}+{\mathcal{B}_2\epsilon {\bf 1}_p})\right)
\end{equation}
\begin{equation}
\Longrightarrow  \left.\frac{\partial J_{R1}}{\partial i^{e\ast}}\right|_{i^{e\ast}=i^{e*(\ell)}}=2{\mathcal{B}_2}^{\top}L{\mathcal{B}_1}.
\end{equation}
Similarly, the gradient of $J_{R2}$  in (\ref{relax}) can be derived as follows. 
Note that for any $x \in \mathbb R^n$, $y \in \mathbb R^n$ and a symmetric matrix $S \in \mathbb R^{n \times n}$, it is known that $\frac{\partial x^{\top}y}{\partial x} = y$ and $\frac{\partial x^{\top}Sx}{\partial x} = 2Sx$. Thus, $\frac{\partial}{\partial x} (x-y)^{\top}S(x-y) = 2S(x-y)$. Therefore, 
\begin{equation}
\left.\frac{\partial J_{R2}}{\partial {i}^{e\ast}}\right|_{{i}^{e\ast}=i^{e*(\ell)}}=2\mathcal D({i}^{e\ast(\ell)}-{i}^{eD}). 
\end{equation}


\bibliographystyle{IEEEtran}
\bibliography{ref}

\begin{thebibliography}{10}
\providecommand{\url}[1]{#1}
\csname url@samestyle\endcsname
\providecommand{\newblock}{\relax}
\providecommand{\bibinfo}[2]{#2}
\providecommand{\BIBentrySTDinterwordspacing}{\spaceskip=0pt\relax}
\providecommand{\BIBentryALTinterwordstretchfactor}{4}
\providecommand{\BIBentryALTinterwordspacing}{\spaceskip=\fontdimen2\font plus
\BIBentryALTinterwordstretchfactor\fontdimen3\font minus \fontdimen4\font\relax}
\providecommand{\BIBforeignlanguage}[2]{{%
\expandafter\ifx\csname l@#1\endcsname\relax
\typeout{** WARNING: IEEEtran.bst: No hyphenation pattern has been}%
\typeout{** loaded for the language `#1'. Using the pattern for}%
\typeout{** the default language instead.}%
\else
\language=\csname l@#1\endcsname
\fi
#2}}
\providecommand{\BIBdecl}{\relax}
\BIBdecl

\bibitem{qian}
K.~Qian, C.~Zhou, M.~Allan, and Y.~Yuan, ``Modeling of {L}oad {D}emand {D}ue to {EV} {B}attery {C}harging in {D}istribution {S}ystems,'' \emph{IEEE Transactions on Power Systems}, vol.~26, no.~2, pp. 802--810, May 2011.

\bibitem{floch}
C.~Le~Floch, F.~Belletti, and S.~Moura, ``{O}ptimal {C}harging of {E}lectric {V}ehicles for {L}oad {S}haping: {A} {D}ual-{S}plitting {F}ramework {w}ith {E}xplicit {C}onvergence {B}ounds,'' \emph{IEEE Transactions on Transportation Electrification}, vol.~2, no.~2, pp. 190--199, June 2016.

\bibitem{yinling}
Y.~Xu, ``{O}ptimal {D}istributed {C}harging {R}ate {C}ontrol of {P}lug-{I}n {E}lectric {V}ehicles for {D}emand {M}anagement,'' \emph{IEEE Transactions on Power Systems}, vol.~30, no.~3, pp. 1536--1545, May 2015.

\bibitem{Bayram}
I.~S. Bayram, G.~Michailidis, M.~Devetsikiotis, and F.~Granelli, ``{E}lectric {P}ower {A}llocation in a {N}etwork of {F}ast {C}harging {S}tations,'' \emph{IEEE Journal on Selected Areas in Communications}, vol.~31, no.~7, pp. 1235--1246, July 2013.

\bibitem{eig_paper_2}
W.~Du, Q.~Fu, and H.~F. Wang, ``{S}mall-{S}ignal {S}tability of a {DC} {N}etwork {P}lanned for {E}lectric {V}ehicle {C}harging,'' \emph{IEEE Transactions on Smart Grid}, vol.~11, no.~5, pp. 3748--3762, September 2020.

\bibitem{eig_paper_3}
Q.~Fu, W.~Du, and H.~Wang, ``Planning of the {DC} {S}ystem {C}onsidering {R}estrictions on the {S}mall-{S}ignal {S}tability of {EV} {C}harging {S}tations and {C}omparison {B}etween {S}eries and {P}arallel {C}onnections,'' \emph{IEEE Transactions on Vehicular Technology}, vol.~69, no.~10, pp. 10\,724--10\,735, July 2020.

\bibitem{new1}
M.~Tabari and A.~Yazdani, ``Stability of a {DC} {D}istribution {S}ystem for {P}ower {S}ystem {I}ntegration of {P}lug-{I}n {H}ybrid {E}lectric {V}ehicles,'' \emph{IEEE Transactions on Smart Grid}, vol.~5, no.~5, pp. 2564--2573, September 2014.

\bibitem{eig_paper_5}
M.~M. Mahfouz and M.~R. Iravani, ``{G}rid-{I}ntegration of {B}attery-{E}nabled {DC} {F}ast {C}harging {S}tation for {E}lectric {V}ehicles,'' \emph{IEEE Transactions on Energy Conversion}, vol.~35, no.~1, pp. 375--385, March 2020.

\bibitem{new2}
Y.~Jin, M.~A. Acquah, M.~Seo, and S.~Han, ``{O}ptimal {S}iting and {S}izing of {EV} {C}harging {S}tation {U}sing {S}tochastic {P}ower {F}low {A}nalysis for {V}oltage {S}tability,'' \emph{IEEE Transactions on Transportation Electrification}, vol.~10, no.~1, pp. 777--794, March 2024.

\bibitem{new3}
C.~Dharmakeerthi, N.~Mithulananthan, and T.~Saha, ``Impact of {E}lectric {V}ehicle {F}ast {C}harging on {P}ower {S}ystem {V}oltage {S}tability,'' \emph{International Journal of Electrical Power \& Energy Systems}, vol.~57, pp. 241--249, May 2014.

\bibitem{evcs_model}
N.~Pogaku, M.~Prodanovic, and T.~C. Green, ``{M}odeling, {A}nalysis and {T}esting of {A}utonomous {O}peration of an {I}nverter-{B}ased {M}icrogrid,'' \emph{IEEE Transactions on Power Electronics}, vol.~22, no.~2, pp. 613--625, March 2007.

\bibitem{LCL_design}
A.~Arancibia and K.~Strunz, ``{M}odeling of an {E}lectric {V}ehicle {C}harging {S}tation for {F}ast {DC} {C}harging,'' in \emph{IEEE International Electric Vehicle Conference}, March 2012, pp. 1--6.

\bibitem{ref1}
S.~Sivakumar, S.~B. Ramasamy~Gunaseelan, M.~V. Reddy~Krishnakumar, N.~Krishnan, A.~Sharma, and A.~Aguila~Téllez, ``{A}nalysis of {D}istribution {S}ystems in the {P}resence of {E}lectric {V}ehicles and {O}ptimal {A}llocation of {D}istributed {G}enerations {C}onsidering {P}ower {L}oss and {V}oltage {S}tability {I}ndex,'' \emph{IET Generation, Transmission \& Distribution}, vol.~18, no.~6, pp. 1114--1132, March 2024.

\bibitem{khalil_book}
H.~K. Khalil, \emph{{N}onlinear {S}ystems}, 3rd~ed.\hskip 1em plus 0.5em minus 0.4em\relax Upper Saddle River, NJ: Prentice Hall, December 2007.

\bibitem{h2pfm}
M.~Inoue, T.~Sadamoto, M.~Arahata, and A.~Chakrabortty, ``{O}ptimal {P}ower {F}low {D}esign for {E}nhancing {D}ynamic {P}erformance: {P}otentials of {R}eactive {P}ower,'' \emph{IEEE Transactions on Smart Grid}, vol.~12, no.~1, pp. 599--611, January 2021.

\bibitem{wait_save}
\BIBentryALTinterwordspacing
L.~Blog. (2020) Wait \& {S}ave: The {M}ost {A}ffordable {Lyft} {R}ide for {H}ouseholds and {I}ndividuals. Accessed on: August 05, 2023. [Online]. Available: \url{https://www.lyft.com/blog/posts/wait-and-save}
\BIBentrySTDinterwordspacing

\bibitem{EPRI}
D.~Symanski and C.~Watkins, ``380{VDC} {D}ata {C}enter at {D}uke {E}nergy,'' 2010.

\bibitem{charging_price}
\BIBentryALTinterwordspacing
M.~Services. (2023) The {A}verage {C}osts of using {C}ar {C}harging {S}tations. Accessed on: September 05, 2023. [Online]. Available: \url{https://www.mach1services.com/costs-of-using-car-charging-stations/}
\BIBentrySTDinterwordspacing

\bibitem{Chakra}
M.~Chakravorty and D.~Das, ``{V}oltage {S}tability {A}nalysis of {R}adial {D}istribution {N}etworks,'' \emph{International Journal of Electrical Power \& Energy Systems}, vol.~23, no.~2, pp. 129--135, February 2001.

\end{thebibliography}
\vspace{-70pt}
 \begin{IEEEbiography}{Amit Kumer Podder}
received the B.S. and M.S. degrees in electrical engineering from the Khulna University of Engineering \& Technology, Khulna, Bangladesh in 2016 and 2019, respectively. He is currently working toward a Ph.D. degree at North Carolina State University in Raleigh, USA. His research interests include power systems dynamic stability analysis and control with high penetration of inverter-based converters. 
\end{IEEEbiography}
\vspace{-70pt}
 \begin{IEEEbiography}{Tomonori Sadamoto}
 received his Ph.D. degree from the Tokyo Institute of Technology, Tokyo, Japan, in 2015. From November 2018 to October 2023, he was an Assistant Professor with the Department of Mechanical and Intelligent Systems Engineering at the University of Electro-Communications. Since November 2023, he has been an Associate Professor with the same department. In 2014, he was named as a finalist of the 13th European Control Conference Best Student-Paper Award. In 2020, he received the IEEE Control Systems Magazine Outstanding Paper Award. 
\end{IEEEbiography}
\vspace{-70pt}
 \begin{IEEEbiography}{Aranya Chakrabortty}
(SM'15) received his PhD in Electrical Engineering from Rensselaer Polytechnic Institute in 2008. He is currently a Professor in the Electrical and Computer Engineering department at North Carolina State University, Raleigh, NC. His research interests are in all branches of control theory with applications to electric power systems. He received the NSF CAREER award in 2011 and was chosen as a Faculty Scholar by the NC State Provost's office in 2019.  
\end{IEEEbiography}
\end{document}